\documentclass[11pt, a4paper, reqno, numberwithinsect,thm-restate,autoref]{amsart}
\usepackage[margin=1.2in]{geometry}
\usepackage{enumerate}

\usepackage{amsthm}

\usepackage{mathtools}

\numberwithin{equation}{section}
\numberwithin{table}{section}

\usepackage[colorlinks=true,linkcolor=blue]{hyperref}

\numberwithin{equation}{section}

\usepackage{amsmath}
\usepackage[mathscr]{euscript}
\usepackage{amssymb}
\usepackage[noadjust]{cite}
\usepackage{comment}
\usepackage{color}
\usepackage{filecontents}
\usepackage{bm}

\makeatletter
\newcommand{\Eqref}[1]{\textup{\tagform@{\ref*{#1}}}}
\makeatother

\DeclareMathAlphabet{\mymathbb}{U}{bbold}{m}{n}

\usepackage[noabbrev]{cleveref}

\newcommand{\myref}[2][]{%
	\hyperref[#2]{\textup{(}\textnormal{\textit{#1}}\textup{)}}%
}

\newtheorem{thm}{Theorem}[section]
\newtheorem{theorem}{Theorem}[section]
\newtheorem{lemma}[thm]{Lemma}

\newtheorem{proposition}[thm]{Proposition}

\newtheorem{corollary}[thm]{Corollary}

\newtheorem{example}[thm]{Example}

\theoremstyle{definition}

\theoremstyle{remark}

\numberwithin{equation}{section}

\newcommand{\smnoind}{\smallskip\noindent}

\newcommand{\ps}{\mathcal{PS}}

\newcommand{\Z}{\mathrm{Z}}

\newcommand{\U}{\mathrm{U}}

\newcommand{\sphe}{\hbox{Sph}_{_{\RS_{\mathfrak{A}^+}}}}

\newcommand{\spheA}{\hbox{Sph}_{_{\RS_{{A}^+}}}}
\newcommand{\spheAq}{\hbox{Sph}_{_{\RS_{U_{q}(\mathfrak{A})^+}}}}

\renewcommand{\emptyset}{\varnothing}
\newcommand{\RB}{\mathrm{B}}
\newcommand{\RS}{\mathrm{S}}
\newcommand{\CA}{\mathfrak{A}}
\newcommand{\CB}{\mathfrak{B}}
\newcommand{\CP}{\mathcal{P}}

\newcommand{\one}{\mathbf{1}}

\numberwithin{equation}{section}

\begin{document}

\title[Identifying JBW$^*$-algebras through their spheres of positive elements]{Identifying JBW$^*$-algebras through their spheres of positive elements}

\author[Peralta, Saavedra]{Antonio M. Peralta \and Pedro Saavedra}

\address[A.M. Peralta]{Instituto de Matem{\'a}ticas de la Universidad de Granada (IMAG), Departamento de An{\'a}lisis Matem{\'a}tico, Facultad de
	Ciencias, Universidad de Granada, 18071 Granada, Spain.}
\email{aperalta@ugr.es}

\address[Pedro Saavedra]{Departamento de An\'{a}lisis Matem\'{a}tico, Facultad de Ciencias, Universidad de Granada, 18071 Granada, Spain.}
\email{psaavedraortiz@ugr.es}



\date{\today}

\begin{abstract} Let $\CA$ and $\CB$ be JBW$^*$-algebras whose  lattices of projections are denoted by $\CP (\CA)$ and $\CP (\CB)$, respectively, and let $\Theta: \CP (\CA)\to \CP(\CB)$ be an order isomorphism. We prove that if $\CA$ does not contain any type $I_2$ direct summand and $\Theta$ preserves pairs of points at diametrical distance (i.e. points at distance $1$), then $\Theta$ extends to a Jordan $^*$-isomorphism from $\CA$ onto $\CB$. We also establish that if $\CA$ and $\CB$ are two atomic JBW$^*$-algebras of type $I_2$ and $\Theta: \CP (\CA)\to \CP(\CB)$ is an order isomorphism preserving points at distance $\frac{\sqrt{2}}{2}$ in both directions, that is, $$
	\hbox{$\|p-q\| = \frac{\sqrt{2}}{2}$ in $\CP(\CA)$ if, and only if, $\|\Theta (p) -\Theta(q) \| = \frac{\sqrt{2}}{2}$ in $\CP(\CB)$,}$$
	then $\CA$ is Jordan $^*$-isomorphic to $\CB$. We exhibit a counterexample showing that preservation of points at distance $1$ is not sufficient in this case. Furthermore, if $\CA$ and $\CB$ are two general JBW$^*$-algebras such that the type $I_2$ part of $\CA$ is atomic and $\Theta$ is an isometry, we prove the existence of a unique extension of $\Theta$ to a Jordan $^*$-isomorphism from $\CA$ onto $\CB$.\smallskip

	We employ these results to provide a positive answer to Tingley's problem for positive spheres showing that if  $\CA$ and $\CB$ are JBW$^*$-algebras such that the type $I_2$ part of $\CA$ is atomic, then every surjective isometry from the set, $\RS_{\mathfrak{A}^+}$, of positive norm-one elements of $\CA$ onto the positive norm-one elements of $\CB$ extends $($uniquely$)$ to a Jordan $^*$-isomorphism from $\CA$ onto $\CB$.\smallskip

	Among the other tools developed in this paper we include a metric characterization of projections in JBW$^*$-algebras in the following sense: if $a$ is a norm-one positive element in a general JBW$^*$-algebra $\mathfrak{A}$, then $a$ is a projection if, and only if, it satisfies the double sphere property, that is,
	$$\Big\{c \in \RS_{\mathfrak{A}^+} : \|c - b\| = 1 \; \text{for all} \; b \in \RS_{\mathfrak{A}^+} \; \text{with} \; \|b - a\| = 1\Big\} = \{a\}.$$
\end{abstract}

\maketitle

\section{Introduction}\label{sec:intro}

The fascinating challenge of determining whether two Banach spaces are isometrically isomorphic (as real Banach spaces) if, and only if, their unit spheres are isometrically isomorphic, has given rise to attractive results on identifying C$^*$-algebras, JB$^*$-algebras, and more generally, JB$^*$-triples. These results are all specific cases of the so-called Tingley's problem (cf. \cite{Tingley,Peralta18, Pe24}). For example, two von Neumann algebras are isometrically isomorphic if, and only if, there exists a surjective isometry between their unit spheres \cite{FP3}. More generally, if $A$ denotes a unital C$^*$-algebra or a real von Neumann algebra, and $X$ is a real or complex Banach space, every surjective isometry from the unit sphere of $A$ onto the unit sphere of $X$ admits an extension to a real linear isometry from $A$ onto $X$ \cite{MO}. The same conclusion holds when $A$ is a JBW$^*$-triple \cite{BeCuFerPe21}, or a unital JB$^*$-algebra \cite{PeSvard23}.  

In the case that $\CA$ is a unital C$^*$-algebra or a JB$^*$-algebra with the unit sphere $\RS_{\CA}$, there are some particular (strictly smaller) subsets of $\RS_{\CA}$ of special attractiveness as metric invariants, such as the set of unitaries, the lattice of all projections, and the positive spheres (i.e., the set of all positive norm-one elements, $\RS_{\CA^+} = \RS_{\CA}\cap {\CA^+}$). Indeed,  every surjective isometry between the sets of unitary elements of two von Neumann algebras (or more generally of two JBW$^*$-algebras) extends to a surjective real linear isometry between the corresponding algebras (see \cite{HatMol2014,CuPe2022,Hat2014,CuEnHiMiPe22}). 

The study of surjective isometries between the projection lattices $\CP(M)$ and $\CP(N)$ of two von Neumann algebras $M$ and $N$, respectively, has been considered in recent studies. A connected component of $\CP(M)$ which contains more than one element is called a \emph{Grassmann space} in $M$. Let $\mathscr{P}\subset M$ and $\mathscr{Q} \subset N$ be proper Grassmann spaces in $M$ and $N$. A result by M. Mori proves that for every surjective isometry $\Delta : \mathscr{P}\to \mathscr{Q}$ there exist a Jordan $^*$-isomorphism $\Phi :M\to N$ and a central projection $r\in \CP(N)$ which satisfy $\Delta (p) = \Phi (p) r + \Phi (\one-p) (\one-r),$ for all $p\in \mathscr{P}$ (see \cite[Theorem 2.1]{mori2019isometries}). The case of $B(H)$ was previously developed in a series of papers by F. Botelho, J. Jamison, and L. Molnár \cite{BJM}, and Gehér and \v{S}emrl \cite{geher2016isometries,geher2018isometries}. A full description of order isomorphisms between effect algebras of atomic JBW$^*$-algebras is given in \cite{RoeWortl2020}.

A third line of study poses the positive spheres of C$^*$-algebras as metric invariants. The main result in \cite{nagy2018isometries} proves that if $H$ is a finite-dimensional complex Hilbert space, every isometry $\Delta : \RS_{B(H)^+}\to \RS_{B(H)^+}$ is surjective and admits an extension to a surjective complex linear isometry on $B(H)$. The first author of this note showed that the same conclusion holds for arbitrary complex Hilbert spaces \cite{Peralta2019unit}. The so-called \emph{Tingley's problem for positive unit spheres} asks whether every surjective isometry between the positive spheres of two C$^*$-algebras, or more generally between two JB$^*$-algebras, extends to a Jordan $^*$-isomorphism between the corresponding algebras \cite{Peralta18}. A complete positive solution to this problem in the case of type $I$ finite von Neumann algebras that have bounded dimensions of irreducible representations is established in \cite[Theorem 4.5]{LN}. C.-W. Leung, C.-K. Ng, and N.-C. Wong recently gave a positive answer to Tingley's problem for positive spheres in the case of two commutative unital C$^*$-algebras (see \cite[Theorem 15]{leung2021variant}). The same authors solved this problem in the case of von Neumann algebras and AW$^*$-algebras \cite{LeungNgWongTAMS}, whereas the case of non-unital commutative C$^*$-algebras was subsequently treated by K. Ezumi, M.-R. Lin, and T. Miura in \cite{EzumiLinMiura2026}. 

In this article we explore the role of the sets $\CP(\CA)\setminus\{0\}$ and $\RS_{\CA^+},$ of all non-zero projections and all positive norm-one elements in a JBW$^*$-algebra $\CA,$ respectively, as metric invariants for $\CA$. We shall show that, besides requiring completely new arguments, the Jordan setting hides enormous differences. Namely, the simple structure of type $I_2$ von Neumann algebras implies that two von Neumann algebras $M$ and $N$ of this type are C$^*$-isomorphic if, and only if, their centres are isomorphic commutative von Neumann algebras, which is also equivalent to the existence of an order isomorphism preserving elements at diametrical distance between their projection lattices (cf. \cite[Proposition 3.2$(d)$]{LeungNgWongTAMS}). Note that the diameter of the set $\RS_{\CA^+}$ is one. The conclusion differs in the case of JBW$^*$-algebras. In \Cref{counterexample preservers of diametrical distance spin} we show the existence of two factor JBW$^*$-algebras of type $I_2$ or spin, $\mathcal{V}_3$ and $\mathcal{V}_4$, which are not Jordan $^*$-isomorphic but we can find an order isomorphism $\Theta: \CP (\mathcal{V}_3) \to \CP (\mathcal{V}_4)$ preserving points at diametrical distance.

In \Cref{Sec:Spheres} we present several characterizations of algebraic notions, like projections and the relations of partial order and orthogonality among them, in terms of distances. The central concept is the unit sphere of positive norm-one elements around a subset $\mathscr{S}$ in $\RS_{\CA^+}$ defined by $$\hbox{Sph}_{_{\RS_{\mathfrak{A}^+}}} (\mathscr{S}) :=\left\{ x\in \RS_{\mathfrak{A}^+} : \|x-s\|=1 \hbox{ for all } s\in \mathscr{S} \right\}.$$ The main result in \cite{peralta2018characterizing} shows how the double sphere around a positive norm-one element can be employed to characterize projections in $B(H)$, and the characterization has been extended to projections in an arbitrary AW$^*$-algebra $A$ in \cite{LeungNgWongTAMS}; concretely, for each element $a\in \RS_{A^+}$ we have \begin{equation}\label{eq chracterization projections AWstar}
	\boxed{ a \hbox{ is a projection in } A}  \Leftrightarrow \boxed{\spheA(\spheA(a)) = \{a\}}.
\end{equation}  One of the central novelties in \Cref{Sec:Spheres} proves that for each positive norm-one element $a$ in a JB$^*$-algebra $\mathfrak{A}$ the equality $\sphe(\sphe(a)) = \{a\}$ implies that $a$ is a projection in $\mathfrak{A}$. In the case that $\mathfrak{A}$ is a JBW$^*$-algebra, the equality $\sphe(\sphe(a)) = \{a\}$ characterizes when $a$ is a projection (see \Cref{prop:lem:metr-discr-proj}).

\Cref{sec:two projections} is entirely devoted to developing a two-projections-theory for JB$^*$-algebras. We present a complete description of the JB$^*$-subalgebra of a JB$^*$-algebra generated by two projections (see \Cref{p JBstar subalgebra generated by two projections}), as well as the JBW$^*$-subalgebra of a JBW$^*$-algebra generated by two projections and the unit element (see \Cref{p JBW generated by two projecitons}). Among the consequences of these {structure results,} it is established in \Cref{lem:characterization-orthogonality JB*-algebra} (respectively, \Cref{lem:describ-disj-proj}) that two non-zero projections $p,q$ in a JB$^*$-algebra (respectively, a JBW$^*$-algebra) $\mathfrak{A}$ are orthogonal if, and only if, $\|a - b\| = 1$ for all $a, b \in S_{\mathfrak{A}^+}$ such that $a = U_p(a)$, $b = U_q(b)$ (respectively, for every $r,s\in \CP(\CA)\setminus \{0\}$ with $r\leq p$ and $s\leq q$ we have $\| r-s\| = 1$). It is also shown that the partial order between projections in a JBW$^*$-algebra $\CA$ can be characterized in terms of the metric space given by the positive sphere of $\CA$ (see \Cref{lem:describ-ord-proj} and \Cref{c: rem:inver}).

\Cref{Sec: preservers of diametrical distance among projections} is devoted to {exploring} the role of the lattice of projections in a JBW$^*$-algebra as a metric invariant. As we have already commented, there exist non-isomorphic JBW$^*$-algebra factors of type $I_2$ whose lattices of projections can be related by an order isomorphism preserving points at diametrical distance (cf. \Cref{counterexample preservers of diametrical distance spin}). However, if $\Theta:\CP(\CA) \to \CP(\CB)$ is an order isomorphism
that preserves points at diametrical distance, that is, $$\hbox{$\|p-q\| =1$ in $\CP(\CA)$ if, and only if, $\|\Theta (p) -\Theta(q) \| =1$ in $\CP(\CB)$,}$$ where $\CA$ and $\CB$ are JBW$^*$-algebras, then $\Theta$ is bi-orthogonality preserving and {maps} central projections in $\CA$ to central projections in $\CB$. If we additionally assume that $\CA$ does not contain any type $I_2$ direct summand, then $\Theta$ extends to a Jordan {$^*$-}isomorphism from $\CA$ onto $\CB$ (see \Cref{prop 4 point 1}). Assuming that $\CA$ and $\CB$ are two atomic JBW$^*$-algebras of type $I_2$, we prove that the existence of an order isomorphism $\Theta: \CP (\CA)\to \CP(\CB)$ preserving points at distance $\frac{\sqrt{2}}{2}$ in both directions, that is, $$\|p-q\| = \frac{\sqrt{2}}{2} \Leftrightarrow \|\Theta (p)-\Theta (q) \| = \frac{\sqrt{2}}{2},$$ implies that $\CA$ and $\CB$ are Jordan $^*$-isomorphic. If the hypothesis on $\Theta$ is replaced with the stronger assumption that $\Theta$ is an isometry, the mapping $\Theta$ extends to a Jordan $^*$-isomorphism from $\CA$ onto $\CB$ (see \Cref{p atomic type I2}). The main conclusion of this section is stated in \Cref{thm:ord-pre-isom}, where we prove that if $\CA$ and $\CB$ are JBW$^*$-algebras, such that the type $I_2$ part of $\CA$ is atomic, and $\Theta:\CP(\CA) \to \CP(\CB)$ is an isometric order isomorphism, then $\Theta$ extends to a Jordan $^*$-isomorphism from $\CA$ onto $\CB$.

In \Cref{Sec:isometries between positive spheres}, our study on Tingley's problem for positive spheres of JBW$^*$-algebras culminates. The conclusion is as follows: Let $\CA$ and $\CB$ be JBW$^*$-algebras such that the type $I_2$ part of $\CA$ is atomic. Then every surjective isometry $\Delta:\RS_{\CA^+} \to \RS_{\CB^+}$ extends $($uniquely$)$ to a Jordan $^*$-isomorphism from $\CA$ onto $\CB$ (see \Cref{thm:met-pre-proj}).

\section{Spheres around sets of positive elements and a metric characterization of projections}\label{Sec:Spheres}


A real (respectively, complex) Jordan-Banach algebra is a real (respectively, complex) Banach space $\CA$ together with a continuous, not necessarily associative, commutative product (denoted by $\circ$) which satisfies the so-called \emph{Jordan identity}: 
\begin{equation}\label{eq Jordan identity}
	(a\circ b)\circ b^2  = (a\circ b^2)\circ b, \ \ (\forall a,b\in \mathfrak{A}).	
\end{equation} Unless otherwise specified, Jordan-Banach algebras in this paper are complex. For each pair of elements $a,c$ in $\mathfrak{A}$, the symbol $U_{a,c}$ will stand for the linear mapping on $\mathfrak{A}$ defined by $U_{a,c}(b) := (a\circ b)\circ c + (b\circ c)\circ a - (a\circ c)\circ b$ ($b\in \mathfrak{A}$). We shall simply write $U_a$ for $U_{a,a}$; equivalently, $U_a(b)=2a\circ(a\circ b)-a^2\circ b$. The mapping $U_a$ is positive whenever $a$ is self-adjoint \cite[Proposition 3.3.6]{HO1983}. In the case that a C$^*$-algebra $A$ is regarded as a JB$^*$-algebra with respect to the natural Jordan product, the mapping $U_a$ is given by $U_a (x) = a x a $ for all $a,x\in A$. The Jordan multiplication operator by a fixed element $a\in \mathfrak{A}$ will be denoted by $M_a$, that is, $M_a (x) = a\circ x$ ($x\in \mathfrak{A}$). If $p$ is a projection, then $U_p$ is the quadratic projection onto the Peirce $2$-subspace $U_p(\mathfrak{A})$; in particular, $U_p(\one)=p$ and $U_p(\mathfrak{A})$ is a JB$^*$-algebra with unit $p$ whenever $\mathfrak{A}$ is a JB$^*$-algebra \cite[\S 2.6]{HO1983}. If $s$ is a symmetry (i.e., $s^2 = \mathbf{1}$, $s^* = s$), then $U_s$ is a surjective Jordan $^*$-automorphism of $\mathfrak{A}$ \cite[2.8.6]{HO1983}.

Every real or complex associative Banach algebra is a real or complex Jordan-Banach algebra, respectively, with respect to the product $a\circ b: = \frac12 (a b +ba)$. There are known examples of Jordan algebras which cannot be embedded as Jordan subalgebras of an associative algebra, for example the Jordan algebra $H_3(\mathbb{O})$ of all Hermitian $3\times 3$ matrices with entries in the algebra of  octonions or Cayley numbers (see \cite[Corollary 2.8.5]{HOS} or \cite[\S 2.5.1]{Cabrera-Rodriguez-vol1} for more information). Jordan algebras of this type are called \emph{exceptional}.

A Jordan-Banach algebra $\mathfrak{A}$ is called \emph{unital} if there exists $\mathbf{1}\in \mathfrak{A}$ satisfying $\one \circ a = a$ for all $a\in \mathfrak{A}$. An element $a$ in a unital Jordan-Banach algebra $\mathfrak{A}$ is called \emph{invertible} whenever there exists $b\in \mathfrak{A}$ satisfying $a \circ b = \one$ and $a^2 \circ b = a.$ The element $b$ is unique and it will be denoted by $a^{-1}$ (cf. \cite[3.2.9]{HOS} and \cite[Definition 4.1.2]{Cabrera-Rodriguez-vol1}). We shall denote by $\mathfrak{A}^{-1}$ the set of all invertible elements in $\mathfrak{A}$. The Jordan-spectrum of an element $a$ in a unital complex Jordan-Banach algebra $\mathfrak{A}$ is the (non-empty compact) set $\hbox{J-}\sigma (a):=\{\lambda\in \mathbb{C}: a-\lambda\one\notin\mathfrak{A}^{-1}\}.$ One of the basic properties of the theory of Jordan algebras, asserts that the closed Jordan subalgebra, $\mathfrak{A}_a$, generated by a single element $a$ (and the unit element) of a complex Jordan-Banach algebra, $\mathfrak{A},$ is an associative and commutative complex Banach algebra (cf. \cite[\S2, Theorem 2.3]{Aupetit95} and  \cite[Theorems 4.1.88 and 4.1.93]{Cabrera-Rodriguez-vol1}). Thus, the usual properties of the set of invertible elements and the holomorphic functional calculus remain valid in the setting of complex Jordan-Banach algebras. 

A \emph{JB$\,^*$-algebra} $\mathfrak{A}$ is a complex Jordan-Banach algebra equipped with an algebra involution $a\mapsto a^*$ satisfying $\left\| U_{a} (a^*) \right\| = \|a\|^3$ ($a\in \mathfrak{A}$). Each C$^*$-algebra $A$ is a JB$^*$-algebra with respect to the original norm and involution and the natural Jordan product. In this particular setting, $U_a (a^* ) = a a^* a$, so the geometric axiom in the definition of JB$^*$-algebras is equivalent to the Gelfand-Naimark axiom. A \emph{JB-algebra} is a real Jordan-Banach algebra $\mathfrak{J}$ satisfying the following axioms: \begin{enumerate}[$(\mathrm{JB}1)$]
	\item $\|a^2\| = \|a\|^2$  for all $a \in \mathfrak{J}$; 
	\item[$(\mathrm{JB}2)$] $\|a^2\| \leq \|a^2 + b^2\|$ for all $a,b \in \mathfrak{J}$.
\end{enumerate}
If $\mathfrak{J}$ is unital, it is clear from $(\mathrm{JB}1)$ that $\|\textbf{1}\| = 1$. 

The classes of JB-algebras and JB$^*$-algebras are mutually determined. {Namely}, for each JB$^*$-algebra $\mathfrak{A}$, the set $\mathfrak{A}_{sa} : = \{a \in \mathfrak{A} \ : \, a^* = a\},$ of all \emph{self-adjoint} or \emph{symmetric} elements in $\mathfrak{A}$, is a JB-algebra \cite[see Proposition 3.8.2]{HOS}, and reciprocally, every JB-algebra corresponds to the self-adjoint part of a (unique) JB$^*$-algebra (see \cite{Wright1977}). Throughout this paper, we shall write $H_3(\mathbb{O}^{\mathbb{C}})$ for the space of all $3\times 3$ Hermitian matrices over the complex Cayley division algebra, that is, the exceptional JB$^*$-algebra whose self-adjoint {part} is $H_3(\mathbb{O})$. This is a finite-dimensional exceptional JBW$^*$-factor of type $I_3$.

A \emph{JBW$\,^*$-algebra} (resp., a \emph{JBW-algebra}) is a JB$^*$-algebra (resp., a JB-algebra) which is also a dual Banach space. Therefore all von Neumann algebras are JBW$^*$-algebras. Every JBW$^*$-algebra is unital (cf. \cite[Lemma 4.1.7]{HOS}).

According to the standard terminology, (norm-closed) JB$^*$-subalgebras of C$^*$-algebras are called \emph{special JB$\,^*$-algebras} or \emph{JC$\,^*$-algebras}. A \emph{JW$\,^*$-algebra} is a JC$^*$-algebra which is also a dual Banach space, or equivalently, a weak$^*$-closed JB$^*$-subalgebra of some von Neumann algebra (see \cite[\S 2]{AlfsenShultz2003}). A milestone result in the theory of JB$^*$-algebras follows from the celebrated Shirshov-Cohn theorem, and affirms that the JB$^*$-subalgebra, $\mathfrak{A}_{a,b},$ of a JB$^*$-algebra, $\mathfrak{A},$ generated by two symmetric elements $a$ and $b$ (and the unit element) is a JC$^*$-algebra, that is, a JB$^*$-subalgebra of some C$^*$-algebra $A$, and in the case that $\mathfrak{A}_{a,b}$ is unital we can {additionally} assume that $\mathfrak{A},$ $\mathfrak{A}_{a,b},$ and {$A$} share the same unit (cf. \cite[Theorem 7.2.5]{HOS} and \cite[Corollary 2.2]{Wright1977}).

The reader interested in further details regarding the theory of JB- and JB$^*$-algebras is {referred} to the monographs \cite{AlfsenShultz2003}, \cite{HOS} and \cite{Cabrera-Rodriguez-vol1}. 

An element $a$ in a JB$^*$-algebra $\mathfrak{A}$ is called \emph{positive} if $a =a^*$ and its Jordan-spectrum is contained in $[0,+\infty)$. The symbol $\mathfrak{A}^{+}$ will denote the cone of all positive elements in $\mathfrak{A}$. Throughout this note we shall write 
$$
\RB_{\mathfrak{A}^+}:= \{a\in \mathfrak{A}^+: \|a\|\leq 1 \} \quad\text{and}\quad
\RS_{\mathfrak{A}^+}:= \{a\in \mathfrak{A}^+: \|a\|= 1 \},
$$
for the \emph{positive unit ball} and the \emph{positive unit sphere} of $\mathfrak{A}$, respectively.
The symbol $\CP({\mathfrak{A}})$ will stand for the set of all projections in $\mathfrak{A}$, that is, the set of all self-adjoint idempotents in $\mathfrak{A}$. To simplify the notation, for each projection $p$ in a unital JB$^*$-algebra $\mathfrak{A}$, we shall write $p^{\perp}$ for $\one-p$. The reader should be warned that we can frequently find JB$^*$-algebras $\mathfrak{A}$ for which $\CP({\mathfrak{A}}) =\{0\}$. 
When $\mathfrak{A}$ is unital, we write $\mathfrak{A}^{-1}$ for the set of all invertible elements in $\mathfrak{A}$, and 
we set
$$
\RB_{\mathfrak{A}^+}^{-1} := \RB_{\mathfrak{A}^+}\cap \mathfrak{A}^{-1}, \quad\text{ and }\quad \RS_{\mathfrak{A}^+}^{-1} := \RS_{\mathfrak{A}^+}\cap \mathfrak{A}^{-1}.
$$ The set of all positive functionals on $\mathfrak{A}$ will be denoted by $\mathfrak{A}^{*}_{+}$. Norm-one positive elements in $\mathfrak{A}^{*}_{+}$ are called \emph{states}. Finally, the set of all \emph{pure states} on $\mathfrak{A}$ (i.e. the extreme points of $\RB_{\mathfrak{A}^{*}_{+}}$) will be denoted by $\ps(A)$. 


\medskip

\begin{lemma}\label{lem:peak-and-zero} Let $\mathfrak{A}$ be a unital JB$\,^*$-algebra,
	$p,q\in \CP (\mathfrak{A}),$ and $a,b\in \RS_{\mathfrak{A}^+}$. Then, the following statements hold:
	\begin{enumerate}[$(a)$]
		\item $\|a-b\| = 1$ if, and only if, there exists $\omega\in \ps(\mathfrak{A})$ such that $\{\omega(a), \omega(b)\} = \{0,1\}$.
		\item $a\in \mathfrak{A}^{-1}$ if, and only if, $\omega(a)\neq 0$ for all $\omega\in\ps(\mathfrak{A})$, which is also equivalent to $\|a - \one \| < 1$. 
		\item  If $\omega\in \ps(\mathfrak{A})$ with $\omega(p)=1$, then $\omega(a\circ p) = \omega (U_p (a)) = \omega(a)$.  
		\item When $\mathfrak{A}$ is unital, one has $\big\{b\in \RS_{\mathfrak{A}^+}^{-1}: p\leq b\big\} = p +  \RB_{\U_{\one-p}(\mathfrak{A})^+}^{-1}$. 	
	\end{enumerate}
\end{lemma}

\begin{proof}
	$(a)$ Suppose that $\|a-b\| =1$. Using Lemma 3.6.8 from \cite{HOS}, we can find $\omega \in \ps(\mathfrak{A})$ such that $1 = \|a - b\| = |\omega(a - b)| = |\omega(a) - \omega(b)|$. Having in mind that $\omega(a), \omega(b) \in [0, 1],$ because $\omega$ is a state, we obtain that $\{\omega(a), \omega(b)\} = \{0, 1\}$. This gives the forward implication. The reciprocal implication is obvious.
	
	$(b)$ Since the JB$^*$-subalgebra of $\mathfrak{A}$ generated by $a$ and $\one$ is a unital commutative C$^*$-algebra in which $a$ is a positive element in the sphere, we deduce that $a$ is invertible if, and only if, $\|\one-a\|<1$, and by statement $(a)$ above (with $b = \one$), the latter is equivalent to $\omega(a)\neq 0$ for all $\omega\in\ps(\mathfrak{A})$.
	
	$(c)$ It is a well-known fact that, if $\omega \in \mathfrak{A}^*$ is any positive functional satisfying $\|\omega\|=\omega(p)$, then $\omega(x) = \omega(x \circ p)$ for all $x \in \mathfrak{A}$ (just observe that by the Cauchy-Schwarz inequality we have $|\omega(x \circ p^{\perp})|^2 \leq \omega (x\circ x^*) \omega (p^{\perp}) = 0$). Even a stronger property holds in the wider setting of JB$^*$-triples (cf. \cite[Proposition 1]{FriedmanRusso1985}).
	
	$(d)$ It is easy to see that $p +  \RB_{\U_{\one-p}(\mathfrak{A})^+}^{-1}\subseteq \big\{b\in \RS_{\mathfrak{A}^+}^{-1}: p\leq b\big\}$.  For the other inclusion, suppose that $b\in \RS_{\mathfrak{A}^+}^{-1}$ with $p\leq b$. Since $b$ is invertible in $\mathfrak{A}$, we know that $\|b - \one\| < 1$, so $\|(b-p)-p^{\perp}\| = \|b - \one\| < 1$. Having in mind $(b)$, we see that $b-p$ is invertible in $\U_{\one-p}(\mathfrak{A})$, and thus $b = p + (b-p) \in p + \RB_{\U_{\one-p}(\mathfrak{A})^+}^{-1}$.
\end{proof}

Recall that elements $a, b$ in a Jordan algebra $\mathfrak{A}$ are said to \emph{operator commute} if $$(a\circ c)\circ b= a\circ (c\circ b),$$ for all $c\in \mathfrak{A}$; equivalently, the mappings $M_a$ and $M_b$ commute in the associative algebra $L(\mathfrak{A})$ of all linear mappings on $\mathfrak{A}$. The \emph{centre} of $\mathfrak{A}$ (denoted by $\Z(\mathfrak{A})$) is defined as the set of all elements in $\mathfrak{A}$ which operator commute with any other element in $\mathfrak{A}$. Elements in $\Z(\mathfrak{A})$ are called central. A classical result by D.M. Topping shows that in the case that $\mathfrak{A}$ is a JC$^*$-algebra, regarded as a JB$^*$-subalgebra of a C$^*$-algebra $A$, two self-adjoint elements $a$ and $b$ in $\mathfrak{A}$ operator commute if, and only if, they commute with respect to the {associative} product of the C$^*$-algebra $A$ (\cite[Proposition 1]{Topping65}). A more recent result proves that if $a$ and $b$ are arbitrary elements in a JB$^*$-algebra $\mathfrak{A}$ satisfying that the JB$^*$-subalgebra of $\mathfrak{A}$ generated by $a$ and $b$ is a JC$^*$-subalgebra of some C$^*$-algebra $A$, then $a$ and $b$ commute in the usual sense as elements of $A$ whenever they operator commute in $\mathfrak{A}$ \cite[Proposition 1.2]{EsPeVill24}.

It is known that $\Z(\mathfrak{A})$ is a commutative C$^*$-algebra. If $\mathfrak{A}$ is a JBW$^*$-algebra, the separate weak$^*$-continuity of the Jordan product of $\mathfrak{A}$ assures that $\Z(\mathfrak{A})$ is weak$^*$-closed in $\mathfrak{A}$, and thus $\Z(\mathfrak{A})$ is a von Neumann algebra. 

Let $\mathfrak{A}$ be a JB$^*$-algebra, and let $\mathscr{S}$ be a subset of $\RS_{\mathfrak{A}^+}$. Following the notation in \cite{peralta2018characterizing,Peralta2019unit}, the \emph{unit sphere around $\mathscr{S}$ in $\RS_{\mathfrak{A}^+}$}\label{def positive sphere around a set} is defined as the set $$\hbox{Sph}_{_{\RS_{\mathfrak{A}^+}}} (\mathscr{S}) :=\left\{ x\in \RS_{\mathfrak{A}^+} : \|x-s\|=1 \hbox{ for all } s\in \mathscr{S} \right\}.$$  In the case that $\mathscr{S}$ reduces to a single element $a$, we shall write $\sphe(a)$ for $\sphe(\{a\})$. 

Let us recall that two projections $p,q$ in a JB$^*$-algebra $\CA$ are called \emph{orthogonal} ($p\perp q$ in short) if $p\circ q=0$. Similarly, two self-adjoint elements in $\CA$ are orthogonal when they have zero product. General elements $a,b\in \CA$ are said to be \emph{orthogonal} ($a\perp b$ in short) if $U_{a,x} (b^*)=0$ for all $x\in \CA$ (see \cite[\S 4]{BurFerGarPe09} for a more detailed exposition). It is known that $a\perp b$ in $\CA$ implies that $\|a + b\| = \max\{\|a\|,\|b\|\}$ (see \cite[Lemma 1.3 $(a)$]{FriedmanRusso1985}).  

In our next lemma we show how some algebraic notions (like operator commutativity of projections and the unit element) can be characterized in {terms} of unit spheres around certain subsets of $\RS_{\mathfrak{A}^+}$.

\begin{lemma}\label{Lemma 2point2} Let $\mathfrak{A}$ be a unital JB$\,^*$-algebra, and let $p$ be a projection in $\mathfrak{A}$. Then the following statements hold:
	\begin{enumerate}[$(a)$]
		\item $p$ is a central projection in $\mathfrak{A}$ if, and only if, $\| p - q\| = 1$ for every $q$ in $\CP (\mathfrak{A})\setminus \{p\}$, equivalently,  $\CP (\mathfrak{A})\setminus \{p\} \subseteq \sphe(p)$. 	
		\item Suppose additionally that $p$ is a central projection in $\mathfrak{A}$. 
		Then $$\hbox{$\boxed{p=\one}$ $\Leftrightarrow$ $\boxed{\hbox{$\sphe (a) \cup \sphe(p) \neq \RS_{\mathfrak{A}^+}$, for every $a\in \RS_{\mathfrak{A}^+}$.}}$} $$
	\end{enumerate}
\end{lemma}

\begin{proof}
	$(a)$ Consider $q\in \CP (\mathfrak{A})\setminus \{p\}$. 
	
	\smnoind
	$\Rightarrow)$ Suppose that $p$ is a central projection in $\mathfrak{A}$. Take any projection $q \in \CP (\mathfrak{A})$ with $q \neq p$. It can be straightforwardly checked from the operator commutativity of $p$ and $q$ that $p \circ (\one - q), q \circ (\one - p) \in \CP (\mathfrak{A})$ with $p \circ (\one - q) \perp q \circ (\one - p)$, and thus
	\begin{equation*}
		\|p-q\| = \|p \circ (\one - q) - q \circ (\one - p)\| = \max\{\|p \circ (\one - q)\|,  \|q \circ (\one - p)\|\}.
	\end{equation*}
	Since $\|p \circ (\one - q)\|,  \|q \circ (\one - p)\| \in \{0, 1\},$ and at least one of these projections is non-zero (otherwise $p = q$, which is impossible), we obtain $\|p-q\| = 1$.
	
	\smnoind
	$\Leftarrow)$ Suppose now that $\|p-q\| = 1$ for every projection $q \in \CP (\mathfrak{A})\setminus \{p\}$. It follows that for each $0 < \delta < 1$ the intersection of the open ball $B(p, \delta)$ of centre $p$ and radius $\delta$ with $\CP (\mathfrak{A})$ reduces to $\{p\}$. Thus, $p$ is norm-isolated in $\CP (\mathfrak{A})$, which is known to be equivalent to $p\in \Z(\mathfrak{A})$ (see, for example, \cite[Proposition 2.2]{CP} and \cite[Proposition 3.1.24]{Cabrera-Rodriguez-vol1}).
	
	\smnoind
	$(b)$ Assume now that $p$ is a central projection. 
	
	$\Rightarrow)$ Suppose that $p = \one$. Put $c:= (\one+a)/2.$
	Since $\| \one -c\| = \|c-a\| = \|(\one-a)/2\| \leq  1/2$, we have $c\notin \sphe(p)$ and $c\notin \sphe (a)$.  
	
	$\Leftarrow)$ Arguing by contradiction, we suppose that $p \neq \one$.
	Since $p$ is central, for any $a \in \RS_{\mathfrak{A}^+}$ we can write 
	\begin{equation*}
		a = a\circ p + a \circ(\one - p) = a_1 + a_2
	\end{equation*}
	with $a_1 = a \circ p \in \U_p(\mathfrak{A})^+$, $a_2 = a \circ (\one - p) \in \U_{\one - p}(\mathfrak{A})^+$, $a_1 \circ a_2 = 0$, and $1 = \|a\| = \max\{\|a_1\|, \|a_2\|\}$. This implies that either $a = a_1 + a_2 \equiv (a_1, a_2) \in \RB_{\U_p(\mathfrak{A})^+} \times \RS_{\U_{\one - p}(\mathfrak{A})^+}$, or $a \equiv (a_1, a_2) \in \RS_{\U_p(\mathfrak{A})^+} \times \RB_{\U_{\one - p}(\mathfrak{A})^+}$. However, 
	$$\RB_{\U_p(\mathfrak{A})^+} \times \RS_{\U_{\one - p}(\mathfrak{A})^+} \subseteq \sphe(p), \quad \text{and} \quad \RS_{\U_p(\mathfrak{A})^+} \times \RB_{\U_{\one - p}(\mathfrak{A})^+} \subseteq \sphe(\one - p).$$ Therefore, $\sphe(p) \cup \sphe(\one - p) = \RS_{\mathfrak{A}^+}$, which is impossible. 
\end{proof}

Let $\mathfrak{A}$ be a JB$^*$-algebra. We shall say that a non-zero projection $p$ in $\mathfrak{A}$ is minimal if $U_p (\mathfrak{A}) = \mathbb{C} p$. When the results on the atomic decomposition in \cite[Proposition 4]{FriedmanRusso1985} particularize to the setting of JB$^*$-algebras, we see that for each pure state 
$\omega\in \ps(A)$ there exists a unique minimal projection $p=p_{\omega}$ in $\mathfrak{A}^{**}$ such that $U_{p_{\omega}} (x) = \omega (x) p$ for all $x\in \mathfrak{A}^{**}$. Furthermore, if we write $\CP_{min} (\mathfrak{A}^{**})$ for the set of all minimal projections in $\mathfrak{A}^{**}$, the mapping $\omega\mapsto p_{\omega}$ is a bijection from $\CP_{min} (\mathfrak{A}^{**})$ onto $\ps(A)$.

Let $p$ be a projection in a JB$^*$-algebra $\mathfrak{A}$, and let $a$ and $b$ be elements in $\mathfrak{A}$ with $a$ positive. It is {known} that $a\perp b $, if and only if, $a\circ b =0$ (cf. \cite[Lemma 4.1]{BurFerGarPe09}). The \emph{inner quadratic annihilator} of a subset $\mathscr{S}\subseteq \mathfrak{A}$ is the set defined by $$^{\perp_q}\mathscr{S} :=\{ a\in M : U_s (a) = 0 \hbox{ for all } s\in \mathscr{S}\}.$$ In the case that $\mathscr{S}\subseteq \mathfrak{A}^{+}$ it is further known (see \cite[Lemma 3.11]{GarLiPeTah2022}) that \begin{equation}\label{eq inner q-annihilator and orthogonality} \hbox{$^{\perp_q}\mathscr{S} \cap \mathfrak{A}^{+} = \mathscr{S}^{\perp} \cap  \mathfrak{A}^{+},$ where $\mathscr{S}^{\perp} =\{x\in \mathfrak{A} : x\perp s \ \forall s\in \mathscr{S}\}.$}
\end{equation}

Another interesting property required for later purposes reads as follows: 
\begin{equation}\label{eq FR Peirce 2 plus orthogonal} \hbox{For all } a\in \RS_{\mathfrak{A}^+}, \ p\in\CP (\mathfrak{A}), \hbox{ the identity } U_p (a) = p \hbox{ implies } a = p + U_{\one-p} (a). 
\end{equation} The statement in \Eqref{eq FR Peirce 2 plus orthogonal} is actually a consequence of \cite[Lemma 1.6]{FriedmanRusso1985}, where a more general conclusion is established.

The set $\CP({\mathfrak{A}})\backslash\{0\},$ of all non-zero projections in a JB$^*$-algebra $\mathfrak{A}${,} is a particular subset of the positive unit sphere of $\mathfrak{A}$.  It is natural to ask whether we can determine the elements in $\CP({\mathfrak{A}})\backslash\{0\}$ as those $a\in \RS_{\mathfrak{A}^+}$ satisfying a metric property. In the case that $A$ is a C$^*$-algebra and $a\in \RS_{{A}^+}$ it is known that  $a$ is a projection in $A$ if $\spheA(\spheA(a)) = \{a\}$ \cite[Proposition 2.2]{peralta2018characterizing}. Recall that for each subset $\mathscr{S}\subseteq  \RS_{\mathfrak{A}^+}$, the set $\sphe(\mathscr{S})$ is a subset of $\RS_{\mathfrak{A}^+}$ defined in {purely} geometric terms (see page \pageref{def positive sphere around a set}). If $A$ is a $B(H)$ space or an atomic von Neumann algebra, Theorem 2.3 in \cite{peralta2018characterizing} assures that $$\boxed{ a\in \RS_{A^+} \hbox{ is a projection in } A}  \Leftrightarrow \boxed{\spheA(\spheA(a)) = \{a\}}.$$ A more recent contribution proves that the same geometric property characterizes non-zero projections in any AW$^*$-algebra $A$ (see \cite[Lemma 4.1]{LeungNgWongTAMS}). The main goal of this section is to extend the mentioned characterization to the case of non-zero projections in a JBW$^*$-algebra.  

\begin{proposition}\label{prop:lem:metr-discr-proj} Let $a$ be a positive norm-one element in a JB$\,^*$-algebra $\mathfrak{A}$.  Then the following statements hold: \begin{enumerate}[$(a)$]
		\item  	$\sphe(\sphe(a)) = \{a\}$ implies that $a$ is a projection in $\mathfrak{A}$.
		\item  If $\mathfrak{A}$ is a JBW$^*$-algebra, the following equivalence holds: $$\boxed{ a \hbox{ is a projection in } \mathfrak{A}}  \Leftrightarrow \boxed{\sphe(\sphe(a)) = \{a\}}.$$ 
	\end{enumerate}
\end{proposition}

\begin{proof} $(a)$ Suppose that $\sphe(\sphe(a)) = \{a\}$. We shall prove that $\hbox{J-}\sigma(a)\subseteq\{0,1\}$. If $\hbox{J-}\sigma(a)\nsubseteq \{0,1\}$, we can find a continuous function $f$ on $[0,1]$ such that $0\leq f\leq 1$, $f(0)=0,$ $f(1) =1,$ and $f(t)\neq t$ for some $t\in \hbox{J-}\sigma(a)$; by continuous functional calculus, $b=f(a)$ is a positive norm-one element with $b \neq a$ (actually $b$ lies in the JB$^*$-subalgebra of $\mathfrak{A}$ generated by $a$). We shall show that such an element $b$ lies in $\sphe(\sphe(a))$, which contradicts the hypothesis.  Having in mind \Cref{lem:peak-and-zero}$(a)$, it suffices to prove that for every $\omega\in \ps(\mathfrak{A})$ with $\omega(a)= 0$ (respectively, $\omega(a)= 1$) we have $\omega(b)= 0$ (respectively, $\omega(b)= 1$). 
	
	Let us take $\omega\in \ps(\mathfrak{A})$ with $\omega(a)= 0$. Let $p_{\omega}$ be the unique minimal projection in $\mathfrak{A}^{**}$ satisfying {$U_{p_{\omega}} (x)= \omega(x) p_{\omega}$ (for all $x\in \mathfrak{A}^{**}$)}. Since, by assumptions $U_{p_{\omega}} (a)=0$, we deduce that $p_{\omega}\perp a$ (cf. \Eqref{eq inner q-annihilator and orthogonality}), and thus $p_{\omega}\perp b$ because $b$ lies in the JB$^*$-subalgebra of $\mathfrak{A}$ generated by $a$. Consequently, $0=U_{p_{\omega}} (b)= \omega(b) p_{\omega}$, and thus $\omega (b) = 0$, as desired. 
	
	Suppose now that $\omega\in \ps(\mathfrak{A})$ with $\omega(a)= 1$, and hence $U_{p_{\omega}} (a)=p_{\omega}$. By applying \Eqref{eq FR Peirce 2 plus orthogonal} we arrive at $a = p_{\omega}+ U_{\one-p_{\omega}} (a),$ with $p_{\omega}\perp  U_{\one-p_{\omega}} (a) $ in $\mathfrak{A}^{**}$. The basic properties of continuous functional calculus show that $b = f(a) =  p_{\omega}+ f(U_{\one-p_{\omega}} (a)),$ with $p_{\omega}\perp  f(U_{\one-p_{\omega}} (a))$ in $\mathfrak{A}^{**}$, and thus $U_{p_{\omega}} (b) = p_{\omega}$ and $\omega(b) =1$.

	\noindent $(b)$ The ``if'' implication is guaranteed by $(a)$. To see the ``only if'' implication, suppose that $p$ is a non-zero projection and
	consider $a\in \sphe(\sphe(p))$. We need to prove that $a=p$. Let us make a simple observation. We consider the positive element $U_p (a)$. If $U_p (a) =0$, we deduce from \Eqref{eq inner q-annihilator and orthogonality} that $a\perp p$, and hence $\|p-a\| =\max\{\|p\|,\|a\|\}= 1$. This proves that $a\in \sphe(p)$, and since $a\in \sphe(\sphe(p))$, we arrive at $0=\|a-a\|=1,$ which is impossible. Therefore $U_p (a)\neq 0$. 
	
	First, we claim that \begin{equation}\label{eq first claim prop charact projections} \| U_q (a)\| =1, \hbox{ for every non-zero projection } q\leq p. 
	\end{equation} Arguing by contradiction, we assume that $\| U_q (a)\| <1$ for some non-zero projection $q\leq p$. The element $U_{\one-q}(a)$  is positive in $\mathfrak{A}$. Observe that $U_{\one-q}(a)=0$ implies that $a\perp \one -q$  (cf. \Eqref{eq inner q-annihilator and orthogonality}), and hence $a = U_q (a)$ {is} norm-one as desired. We can therefore assume that $0<\delta = \|U_{\one-q}(a)\|\leq 1$.  Let $\mathfrak{A}_{U_{\one-q}(a)}$ denote the JB$^*$-subalgebra of $\mathfrak{A}$ generated by $U_{\one-q}(a)$. It is known that $\mathfrak{A}_{U_{\one-q}(a)}$ is JB$^*$-isomorphic to $C_0(\hbox{J-}\sigma(U_{\one-q}(a)))$. Take $b=\frac{1}{\delta}U_{\one-q}(a)$, $x= \frac{\one-q + b}{2}\in \RS_{\mathfrak{A}^+},$ and observe that by $\one-q \geq x, \one- p$, an orthogonality argument gives $$\| p-x \|= \left\| q+ (p-q) -x \right\| = \max\{\|q\|,\|(p-q) -x\|\}=1,$$ that is, $x\in \sphe(p)$, and hence $\| a - x\| =1$. By Lemma \ref{lem:peak-and-zero}(a), there exists $\omega\in \ps(\CA)$ satisfying
	$$
	\{\omega(a), \omega(x)\}=\{0,1\}.
	$$
	If $\omega(x) = 1$, then it follows from $\omega(\one-q) + \omega(b)  = 2$ that $\omega(b) = 1 =\omega(\one-q)$ (observe that $\omega(\one-q), \omega(b)\in [0,1]$), and consequently, the minimal projection $p_{\omega}\in \mathfrak{A}^{**}$ satisfies $p_{\omega}\leq \one-q$ and $U_{p_{\omega}} (b) = p_{\omega}$. It follows from \Eqref{eq FR Peirce 2 plus orthogonal} that $b = p_{\omega} + U_{\one-p_{\omega}} (b)$ with $p_{\omega} \perp U_{\one-p_{\omega}} (b)$ in $\mathfrak{A}^{**}$. Multiplying the identity $b = p_{\omega} + U_{\one-p_{\omega}} (b)$ by $\delta$, we get $$U_{\one-q}(a)=\delta b = \delta p_{\omega}+\delta U_{\one-p_{\omega}} (b),$$ with $\delta U_{\one-p_{\omega}} (b)\in U_{\one-p_{\omega}}(\mathfrak{A}^{**})$. We therefore deduce that $0=\omega (a) = \omega  (U_{\one-q}(a)) = \delta >0$, which is impossible.
	
	If $\omega(x) = 0$,  we necessarily have $\omega(\one-q) = \omega(b)  =0,$ and thus $\omega (q) =1$ and $p_{\omega} \leq q\leq p$. This gives 
	$\omega(a) = \omega(U_q (a)) \leq \|U_q (a)\| < 1,$ contradicting that $\omega (a) =1$. This concludes the proof of the first claim.

	Our next goal will consist in proving that $U_p (a)$ is a projection. Note that as a consequence of the first claim $U_p (a)\in \RS_{U_p(\mathfrak{A})^+}\subseteq \RS_{\mathfrak{A}^+}$. Let $\mathfrak{A}_{U_{p}(a)}$ denote the unital JBW$^*$-subalgebra of $\mathfrak{A}$ generated by $U_{p}(a)$ and $p$. It is known that $\mathfrak{A}_{U_{p}(a)}$ is a commutative von Neumann algebra. If $U_p (a)$ is not a projection, we can find a non-zero projection $q\in \mathfrak{A}_{U_{p}(a)}$ satisfying $\|U_q (a)\| = \|U_q ({U_{p}(a)})\| <1$ (the first equality holds because we are working in the associative von Neumann algebra $\mathfrak{A}_{U_{p}(a)}$), which contradicts the conclusion in the first claim since necessarily ${q}\leq p.$ Therefore, $U_p (a)$ is a projection, and clearly $U_p (a)\leq p$.
	
	If $r= p-U_p (a)\neq 0$, we would have $0 = U_r U_p (a) = U_r (a)$, contradicting again the conclusion in the first claim. This shows that $U_p (a) =p$, and \Eqref{eq FR Peirce 2 plus orthogonal} assures that $a = p + U_{\one-p} (a)$. 
	
	We shall finally show that $U_{\one-p} (a) =0$. Otherwise, the element $x =p +\frac{1}{\| U_{\one-p} (a)\|} U_{\one-p} (a)$ lies in $ \RS_{\mathfrak{A}^+}$ and satisfies $\|p -x\| = 1$, but $$\| a - x\| = \left\| U_{\one-p} (a) - \frac{1}{\| U_{\one-p} (a)\|} U_{\one-p} (a)\right\| = 1- \| U_{\one-p} (a)\| <1,$$ which is impossible. 
\end{proof}

The next corollary can be now obtained as a straightforward consequence of our previous results.    

\begin{corollary}\label{c a surjective isometry maps projections to projections} Let $\CA$ and $\CB$ be JBW$\,^*$-algebras, and let $\Delta:\RS_{\CA^+} \to \RS_{\CB^+}$ be a surjective isometry. Then, $\Delta|_{\CP (\CA)}: \CP (\CA) \to \CP (\CB)$ is a surjective isometry mapping central projections in $\CA$ onto central projections in $\CB$. 
\end{corollary}

\begin{proof} By observing, once again, that the characterization of projections established in the previous \Cref{prop:lem:metr-discr-proj} is entirely determined by the elements in the positive unit sphere and a condition given in terms of distances, we deduce that $\Delta (\CP (\CA)) = \CP (\CB)$. Clearly $\Delta|_{\CP (\CA)}: \CP (\CA) \to \CP (\CB)$ is a surjective isometry. \Cref{Lemma 2point2}$(a)$ gives the final statement. 
\end{proof}

\section{Two projections theory for JB$^*$-algebras}\label{sec:two projections}

We currently have a vast literature devoted to the study of different associative algebras of operators generated by a pair of projections. G.K. Pedersen established the first description of the C$^*$-algebra generated by two projections (see \cite[\S 3]{Ped68}, and the subsequent rediscovery in \cite{VasSpit81,KruRothSilb}). R. Giles and H. Kummer studied the W$^*$-algebra generated by two projections in \cite{GilesKummer71} (see also \cite[Exercise 12.4.11]{KadRingVol2} and \cite[pages 306-308]{Tak}). There are over a {hundred} references on this topic (see for example the excellent survey \cite{BottSpit2010}, the book \cite{RochSantosBook}, and the references therein). In 2015, J. Hamhalter extended the study to the AW$^*$-subalgebra generated by two projections (cf. \cite[Proposition 2.5]{Ham}). Until now no attention has been paid to the JB$^*$-subalgebra generated by two projections. We shall complete the study in this section.

We begin with a result representing the JB$^*$-algebra generated by two projections. Let us introduce some notation. 
Throughout this note we shall write $M_2 (\mathbb{C})$ for the space of all $2\times 2$ complex matrices, and $S_2(\mathbb{C})$ for its subspace of all $2\times 2$ symmetric complex matrices. For each compact Hausdorff space $K$ and each pair of closed subsets $M_1,M_2\subseteq K$ we shall write $C_{_{M_1}}^{_{M_2}}(K, M_2 (\mathbb{C}))$ (respectively, $C_{_{M_1}}^{_{M_2}}(K, S_2 (\mathbb{C}))$) for the C$^*$-algebra of all continuous functions $\mathbf{a} : K\to M_2 (\mathbb{C})$ (respectively, $\mathbf{a} : K\to S_2 (\mathbb{C})$) such that $\mathbf{a}(m)$ is a diagonal matrix for all $m\in M_1$, and $\mathbf{a}_{22} (m) = 0$ for all $m\in M_2$. We shall simply write $C_{_{M_1}}(K, M_2 (\mathbb{C}))$ and  $C_{_{M_1}}(K, S_2 (\mathbb{C}))$ for $C_{_{M_1}}^{\emptyset}(K, M_2 (\mathbb{C}))$ and $C_{_{M_1}}^{\emptyset}(K, S_2 (\mathbb{C})),$ respectively. 

\begin{proposition}\label{p JBstar subalgebra generated by two projections} Let $p$ and $q$ be projections in a JB$\,^*$-algebra  $\mathfrak{A}$, and let $\mathfrak{B}$ denote the JB$\,^*$-subalgebra of $\mathfrak{A}$ generated by $p$ and $q$.  Then for each $t\in \hbox{J-}\sigma(p+2q)\cap \{1, 2, 3\},$ there is a Jordan $^*$-homomorphism $\varphi_t : \mathfrak{B}\to \mathbb{C}$ satisfying $\varphi_1 (p)=1$, $\varphi_1(q)=0,$ $\varphi_2(q)=1,$ $\varphi_2(p)=0,$ and $\varphi_3(p)=\varphi_3(q)=1$, and for each $t\in \hbox{J-}\sigma(U_p (q))\cap (0,1)$ there exists a Jordan $^*$-homomorphism $\pi_t : \mathfrak{B}\to S_2(\mathbb{C})$ such that $\pi_t (p) = \left(\begin{matrix} 1 & 0 \\
		0 & 0
	\end{matrix}\right)$, and $\pi_t (q) = \left(\begin{matrix} t & \sqrt{t (1-t)} \\
		\sqrt{t (1-t)} & 1-t
	\end{matrix}\right)$. Furthermore, one of the following statements holds:
	\begin{enumerate}[$(a)$]
		\item If $0$ and $1$ are not accumulation points of $\hbox{J-}\sigma(U_p(q))$, the  
		JB$\,^*$-algebra $\mathfrak{B}$ is Jordan $^*$-isomorphic to
		$C(\hbox{J-}\sigma(U_p(q))\backslash\{0, 1\},S_2(\mathbb{C}))\oplus^{\infty} C(\hbox{J-}\sigma(p+2q)\cap\{1, 2, 3\})$.
		\item If $0$ is an accumulation point of $\hbox{J-}\sigma(U_p(q))$ but $1$ is not, the  
		JB$\,^*$-algebra $\mathfrak{B}$ is Jordan $^*$-isomorphic to $C_{_{\{0\}}}(\hbox{J-}\sigma(U_p(q))\backslash\{ 1\},S_2(\mathbb{C}))\oplus^{\infty} C(\hbox{J-}\sigma(p+2q)\cap\{3\})$.
		\item If $1$ is an accumulation point of $\hbox{J-}\sigma(U_p(q))$ but $0$ is not, the  
		JB$\,^*$-algebra $\mathfrak{B}$ is Jordan $^*$-isomorphic to $C_{_{\{1\}}}^{_{\{1\}}}(\hbox{J-}\sigma(U_p(q))\backslash\{ 0\},S_2(\mathbb{C}))\oplus^{\infty} C(\hbox{J-}\sigma(p+2q)\cap\{1,2\})$.
		\item If $0$ and $1$ are accumulation points of $\hbox{J-}\sigma(U_p(q))$, the  
		JB$\,^*$-algebra $\mathfrak{B}$ is Jordan $^*$-isomorphic to $C_{_{\{0,1\}}}^{_{\{1\}}}(\hbox{J-}\sigma(U_p(q)),S_2(\mathbb{C}))$.
	\end{enumerate} In all the cases, the Jordan $^*$-isomorphism is given
	by the corresponding restriction of the mapping $\Phi(a)=(\Phi_1(a), \Phi_2(a)),$ where 
	\begin{equation}\label{eq def of Phi1}
		\Phi_1 (a) (t) :=\begin{cases} \pi_t (a), & \hbox{ if } t\in \hbox{J-}\sigma(U_p(q))\cap (0,1),\\
			\left(\begin{matrix} \varphi_1(a) & 0 \\
				0 & \varphi_2 (a)
			\end{matrix}\right), & \hbox{ if }  t=0,\\
			\left(\begin{matrix} \varphi_3(a) & 0\\
				0 & 0
			\end{matrix}\right), & \hbox{ if }  t =1,
		\end{cases}
	\end{equation} and $\Phi_2(a)(t)=\varphi_t (a)$ if $t\in \hbox{J-}\sigma(p+2q)\cap\{1, 2, 3\}$. Furthermore, in the case that $\mathfrak{A}$ is unital and $\mathfrak{B}_{\one}$ denotes the JB$\,^*$-subalgebra of $\mathfrak{A}$ generated by $\{p,q,\one\},$ $\pi_t (\one) = \left(\begin{matrix} 1 & 0 \\
		0 & 1
	\end{matrix}\right)$ for all $t\in \hbox{J-}\sigma(U_p(q))\cap (0,1),$ and there exists an additional Jordan $^*$-homomorphism $\varphi_0 : \mathfrak{B}_{\one}\to \mathbb{C}$ satisfying $\varphi_0 (p)=\varphi_0(q)=0,$ $\varphi_0(\one)=1,$ and one of the following statements holds: 
	\begin{enumerate}[$(a.1)$]
		\item[$(a.1)$] If $0$ and $1$ are not accumulation points of $\hbox{J-}\sigma(U_p(q))$, the  
		JB$\,^*$-algebra $\mathfrak{B}_{\one}$ is Jordan $^*$-isomorphic to
		$C(\hbox{J-}\sigma(U_p(q))\backslash\{0, 1\},S_2(\mathbb{C}))\bigoplus^{\infty} C(\hbox{J-}\sigma(p+2q)\cap\{0, 1, 2, 3\})$.
		\item[$(b.1)$] If $0$ is an accumulation point of $\hbox{J-}\sigma(U_p(q))$ but $1$ is not, the  
		JB$\,^*$-algebra $\mathfrak{B}_{\one}$ is Jordan $^*$-isomorphic to $C_{_{\{0\}}}(\hbox{J-}\sigma(U_p(q))\backslash\{ 1\},S_2(\mathbb{C}))\bigoplus^{\infty} C(\hbox{J-}\sigma(p+2q)\cap\{0,3\})$.
		\item[$(c.1)$] If $1$ is an accumulation point of $\hbox{J-}\sigma(U_p(q))$ but $0$ is not, the  
		JB$\,^*$-algebra $\mathfrak{B}_{\one}$ is Jordan $^*$-isomorphic to $C_{_{\{1\}}} (\hbox{J-}\sigma(U_p(q))\backslash\{ 0\},S_2(\mathbb{C}))\bigoplus^{\infty} C(\hbox{J-}\sigma(p+2q)\cap\{1,2\})$.
		\item[$(d.1)$] If $0$ and $1$ are accumulation points of $\hbox{J-}\sigma(U_p(q))$, the  
		JB$\,^*$-algebra $\mathfrak{B}_{\one}$ is Jordan $^*$-isomorphic to $C_{_{\{0,1\}}}(\hbox{J-}\sigma(U_p(q)),S_2(\mathbb{C}))$,
	\end{enumerate} where in the unital case, the Jordan $^*$-isomorphism is given
	by the corresponding restriction of the mapping $\Phi(a)=(\Phi_1(a), \Phi_2(a)),$ where 
	$$\Phi_1 (a) (t) :=\begin{cases} \pi_{t} (a), & \hbox{ if } t\in \hbox{J-}\sigma(U_p(q))\cap (0,1),\\
		\left(\begin{matrix} \varphi_1(a) & 0 \\
			0 & \varphi_2 (a)
		\end{matrix}\right), & \hbox{ if }  t=0,\\
		\left(\begin{matrix} \varphi_3(a) & 0\\
			0 & \varphi_0(a)
		\end{matrix}\right), & \hbox{ if }  t =1,
	\end{cases}$$ and $\Phi_2(a)(t)=\varphi_t(a)$ for $t\in \hbox{J-}\sigma(p+2q)\cap\{0,1,2,3\}$.
\end{proposition}	

\begin{proof} We shall only prove the non-unital case; the other case follows via similar arguments. By the Shirshov-Cohn theorem \cite[Theorem 1.14]{AlfsenShultz2003} (see also \cite[Theorem 7.2.5]{HOS} or \cite[Corollary 2.2]{Wright1977}), the JB$^*$-algebra $\mathfrak{B}$ is a JC$^*$-algebra, that is, there exists a C$^*$-algebra $A$ containing $\mathfrak{B}$ as a JB$^*$-subalgebra. The elements $p$ and $q$ can be thus regarded as projections in $A$. Let $B$ denote the C$^*$-subalgebra of $A$ generated by $p$ and $q$. The product in $A$ (and $B$) will be simply denoted by mere juxtaposition. According to these assumptions, we have $U_p (q) = p q p$.  Clearly, $\mathfrak{B}$ is a JB$^*$-subalgebra of $B$. It is well-known from spectral theory that $\hbox{J-}\sigma_{\mathfrak{A}}(U_p(q)) = \hbox{J-}\sigma_{\mathfrak{B}}(U_p(q)) = \hbox{J-}\sigma_{B}(U_p(q)) = \hbox{J-}\sigma_{A}(U_p(q)) = \hbox{J-}\sigma_{A}(p q p ) = \sigma_{A}(p q p )\subseteq [0,1]$ (\cite[comments before Lemma 1.23]{AlfsenShultz2003} and \cite[Theorem 4.1.71 and comments before Definition 4.1.2]{Cabrera-Rodriguez-vol1}). Similarly $\hbox{J-}\sigma_{\mathfrak{A}}(p+2q)=\hbox{J-}\sigma_{\mathfrak{B}}(p+2q)=\sigma_{{A}}(p+2q)=\sigma_{B}(p+2q).$ 
	
	By \cite[Theorem 5b]{KruRothSilb} (alternatively, \cite[Theorem 3.4]{Ped68}, \cite[Theorem 4.6]{BottSpit2010}), for each $t\in \hbox{J-}\sigma(p+2q)\cap \{1, 2, 3\},$ there is a $^*$-homomorphism $\psi_t : {B}\to \mathbb{C}$ satisfying $\psi_1 (p)=1$, $\psi_1(q)=0,$ $\psi_2(q)=1,$ $\psi_2(p)=0,$ and $\psi_3(p)=\psi_3(q)=1$, and for each $t\in \hbox{J-}\sigma(U_p (q))\cap (0,1)$ there exists a $^*$-homomorphism $\pi_t : {B}\to M_2(\mathbb{C})$ such that $\pi_t (p) = \left(\begin{matrix} 1 & 0 \\
		0 & 0
	\end{matrix}\right)$, and $\pi_t (q) = \left(\begin{matrix} t & \sqrt{t (1-t)} \\
		\sqrt{t (1-t)} & 1-t
	\end{matrix}\right)$. Moreover, one of the following statements holds:
	\begin{enumerate}[$(1)$]
		\item If $0$ and $1$ are not accumulation points of $\sigma(p q p )= \sigma(U_p(q))$, the  
		C$^*$-algebra ${B}$ is $^*$-isomorphic to
		$C(\sigma(U_p(q))\backslash\{0, 1\},M_2(\mathbb{C}))\bigoplus^{\infty} C(\sigma(p+2q)\cap\{1, 2, 3\})$.
		\item If $0$ is an accumulation point of $\sigma(U_p(q))$ but $1$ is not, the  
		C$^*$-algebra ${B}$ is $^*$-isomorphic to $C_{_{\{0\}}}(\sigma(U_p(q))\backslash\{ 1\},M_2(\mathbb{C}))\bigoplus^{\infty} C(\sigma(p+2q)\cap\{3\})$.
		\item If $1$ is an accumulation point of $\sigma(U_p(q))$ but $0$ is not, the  
		C$^*$-algebra ${B}$ is $^*$-isomorphic to $C_{_{\{1\}}}^{_{\{1\}}}(\sigma(U_p(q))\backslash\{ 0\},M_2(\mathbb{C}))\bigoplus^{\infty} C(\sigma(p+2q)\cap\{1,2\})$.
		\item In the case that $0$ and $1$ are accumulation points of $\sigma(U_p(q))$, the  
		C$^*$-algebra ${B}$ is $^*$-isomorphic to $C_{_{\{0,1\}}}^{_{\{1\}}}(\sigma(U_p(q)),M_2(\mathbb{C}))$.
	\end{enumerate} In all the cases, the $^*$-isomorphism is given
	by the corresponding restriction of the mapping $\Psi(a)=(\Psi_1(a), \Psi_2(a)),$ where 
	$$\Psi_1 (a) (t) :=\begin{cases} \pi_t (a), & \hbox{ if } t\in \sigma(U_p(q))\cap (0,1),\\
		\left(\begin{matrix} \psi_1(a) & 0 \\
			0 & \psi_2 (a)
		\end{matrix}\right), & \hbox{ if }  t=0,\\
		\left(\begin{matrix} \psi_3(a) & 0\\
			0 & 0
		\end{matrix}\right), & \hbox{ if }  t =1,
	\end{cases}$$ and $\Psi_2(a)(t)=\psi_t (a)$ if $t\in \sigma(p+2q)\cap\{1, 2, 3\}$.
	
	Clearly $\Psi (\mathfrak{B})$ coincides with the JB$^*$-subalgebra of $\Psi (B)$ generated by $\Psi(p)$ and $\Psi(q)$, and for each $t\in \hbox{J-}\sigma(p+2q)\cap \{1, 2, 3\},$ the mapping $\varphi_t = \psi_t|_{\mathfrak{B}}$ is a Jordan $^*$-homomorphism, while for each $t\in \hbox{J-}\sigma(U_p(q))\cap (0,1)$ the map $\pi_t|_{\mathfrak{B}} : {\mathfrak{B}}\to M_2(\mathbb{C})$ is a Jordan $^*$-homomorphism too. We shall show next that the cases $(1),$ $(2),$ $(3),$ and $(4)$ lead to $(a),$ $(b),$ $(c),$ and $(d)$, respectively. We shall only detail a couple of cases, for example, the third and the first.
	
	\noindent $(3)$ Suppose $1$ is an accumulation point of $\hbox{J-}\sigma(U_p(q))$ but $0$ is not, and $$\Psi : {B}\to C_{_{\{1\}}}^{_{\{1\}}}(\hbox{J-}\sigma(U_p(q))\backslash\{ 0\},M_2(\mathbb{C}))\oplus^{\infty} C(\hbox{J-}\sigma(p+2q)\cap\{1,2\})$$ is a $^*$-isomorphism. By observing that $C_{_{\{1\}}}^{_{\{1\}}}(\hbox{J-}\sigma(U_p(q))\backslash\{ 0\},S_2(\mathbb{C}))\oplus^{\infty} C(\hbox{J-}\sigma(p+2q)\cap\{1,2\})$ is a JB$^*$-subalgebra of $\Psi(B)$ containing $\Psi (p) = \tilde{p}$ and $\Psi (q)= \tilde{q}$, we deduce that $$\Psi (\mathfrak{B}) \subseteq C_{_{\{1\}}}^{_{\{1\}}}(\hbox{J-}\sigma(U_p(q))\backslash\{ 0\},S_2(\mathbb{C}))\oplus^{\infty} C(\hbox{J-}\sigma(p+2q)\cap\{1,2\})= \tilde{\mathfrak{B}},$$ and consequently 
	$\Psi (\mathfrak{B}) (t)\in S_2 (\mathbb{C})$ for all $t\in \hbox{J-}\sigma(U_p(q))\backslash\{ 0\}$. 
	So, in order to prove $(c)$ it suffices to prove that $\Psi (\mathfrak{B}) = \tilde{\mathfrak{B}}$.
	
	Let $\tilde{r} (t)= U_{\tilde{p}}(\tilde{q}) (t) = \begin{pmatrix}
		t & 0 \\
		0  & 0
	\end{pmatrix} \in \Psi_1 (\mathfrak{B}) \subseteq  C_{_{\{1\}}}^{_{\{1\}}}(\hbox{J-}\sigma(U_p(q))\backslash\{ 0\},S_2(\mathbb{C}))$, and let
	\begin{equation*}
		\tilde{v}_n (t) = \begin{pmatrix}
			0 & t^n \sqrt{t(1-t)} \\
			t^n \sqrt{t(1-t)} & 0
		\end{pmatrix} \in C_{_{\{1\}}}^{_{\{1\}}}(\hbox{J-}\sigma(U_p(q))\backslash\{ 0\},S_2(\mathbb{C})) \quad (n \in \mathbb{N}).
	\end{equation*}
	Then $2 (\tilde{q} - \tilde{r}) \circ \tilde{r} = \tilde{v}_1 \in \Psi_1 (\mathfrak{B})$ and $2 \tilde{v}_n \circ \tilde{r} = \tilde{v}_{n+1}$ for all $n \geq 1$. A simple induction argument allows us to deduce that $\tilde{v}_n \in \Psi_1 (\mathfrak{B})$ for all $n \geq 1$.
	
	Let us consider the closed subspace of functions
	\begin{equation*}
		Z = \left\{\mathbf{a}_{12} (t) : \mathbf{a} = \begin{pmatrix}
			\mathbf{a}_{11} & \mathbf{a}_{12} \\
			\mathbf{a}_{21} & \mathbf{a}_{22}
		\end{pmatrix} \in \Psi_1 (\mathfrak{B}), \mathbf{a}_{11}(t) = \mathbf{a}_{22}(t) = 0 \text{ for all } t \in \hbox{J-}\sigma(U_p(q))\backslash\{ 0\} \right\}.
	\end{equation*} Clearly, $Z \subseteq C_{0}(\hbox{J-}\sigma(U_p(q))\backslash\{ 0,1\})$, and $\overline{z}\in Z$ for all $z\in Z$, since $\Psi_1 (\mathfrak{B})$ is self-adjoint. We will check that $Z = C_{0}(\hbox{J-}\sigma(U_p(q))\backslash\{ 0,1\}) :=\{f\in C(\hbox{J-}\sigma(U_p(q))\backslash\{ 0\}) : f (1) =0\}$.  Observe first that, since $\tilde{v}_n \in \Psi_1 (\mathfrak{B})$ for all $n \geq 1$, it follows that $h(t):=t^n (t - t^2)^{\frac{1}{2}}$ ($t\in \hbox{J-}\sigma(U_p(q))\backslash\{ 0\}$) is a function in $Z$ for all natural $n$. Furthermore, since $\Psi (\mathfrak{B})$ is a JB$^*$-algebra, it can be easily checked that for every $g, h, k \in Z$, the functions $\bar{g}$ and $g \bar{h} k $ also lie in $Z$. In other words, $Z$ is a JB$^*$-subtriple of $C_{0}(\hbox{J-}\sigma(U_p(q))\backslash\{ 0,1\})$. Consequently, for each $z \in Z$ and every odd polynomial  $p(\zeta)$, the element $p(z) (t) = p(z(t)) \in Z$. An application of the Stone-Weierstrass theorem assures that $\tilde{u}_{k} (t) = \sqrt[2k - 1]{h(t) } \in Z$ for all natural $k$. Clearly, $(\tilde{u}_{k})_k$ converges increasingly, pointwise and uniformly to the characteristic function of each compact subset of $\hbox{J-}\sigma(U_p(q))\cap (0,1)$. Having in mind that $Z \subseteq C_{0}(\hbox{J-}\sigma(U_p(q))\backslash\{ 0,1\})$, it is routine to check that 
	\begin{equation*}
		\|g \overline{\tilde{u}_k} h - g h\|_\infty \to 0 \text{ for all } g, h \in Z.
	\end{equation*}  Since $Z$ is a closed subspace, and $g \overline{\tilde{u}_k} h\in Z$ for every natural $k$, the function $gh $ lies in $Z$ for all $g, h \in Z$, which reveals that $Z$ is actually a C$^*$-subalgebra of $ C_{0}(\hbox{J-}\sigma(U_p(q))\backslash\{ 0,1\})$. 
	
	Now, by applying that the family $\{k_n(t) = t^n \sqrt{t - t^2} \ \ (t\in \hbox{J-}\sigma(U_p(q))\backslash\{ 0\})\ : n\in \mathbb{N}\} \subseteq Z$ separates the points in $\hbox{J-}\sigma(U_p(q))\backslash\{ 0,1\}$ and vanishes nowhere, the Stone-Weierstrass theorem leads to $Z = C_{0}(\hbox{J-}\sigma(U_p(q))\backslash\{ 0,1\})$. In other words, for any function $g \in C_{0}(\hbox{J-}\sigma(U_p(q))\backslash\{ 0,1\})$, the element $\begin{pmatrix}
		0 & g(t) \\
		g(t) & 0
	\end{pmatrix}$ lies in $\Psi_1 (\mathfrak{B})$.
	On the other hand, since the elements $U_{\tilde{p}}(\tilde{q}) (t) = \begin{pmatrix}
		t & 0 \\
		0  & 0
	\end{pmatrix},$ $\tilde{q} (t) =  \begin{pmatrix}
		t & \sqrt{t(1-t)} \\ \sqrt{t(1-t)} & 1-t
	\end{pmatrix},$ and $(2 \tilde{p}\circ \tilde{q} -2 U_{\tilde{p}}(\tilde{q})) (t) =  \begin{pmatrix}
		0 & \sqrt{t(1-t)} \\ \sqrt{t(1-t)} & 0
	\end{pmatrix}$ are all in $\Psi_1 (\mathfrak{B})$, we can easily deduce that $$ U_{\tilde{p}}(\tilde{q})^n (t) = \begin{pmatrix}
		t^n & 0 \\ 0 & 0
	\end{pmatrix} , \tilde{x}_n (t) = \begin{pmatrix}
		0 & 0 \\ 0 & (1 - t)^n
	\end{pmatrix} \in \Psi_1 (\mathfrak{B}) \quad \text{for all } n \in \mathbb{N}.$$
	A new application of the Stone-Weierstrass Theorem yields that $\begin{pmatrix}
		f(t) & 0 \\ 0 & 0
	\end{pmatrix}, \begin{pmatrix}
		0 & 0 \\ 0 & h(t)
	\end{pmatrix}$ lie in $\Psi_1 (\mathfrak{B})$ for any $f \in C(\hbox{J-}\sigma(U_p(q))\backslash\{0\})$, $h \in C_{0}(\hbox{J-}\sigma(U_p(q))\backslash\{ 0,1\})$. From this, we can immediately conclude that $\tilde{\mathfrak{B}} \subseteq \Psi (\mathfrak{B})$, as we claimed. The desired Jordan $^*$-isomorphism $\Phi = (\Phi_1,\Phi_2)$ is just the restriction of $\Psi$ to $\mathfrak{B}$.   
	
	\noindent $(1)$ Suppose $0$ and $1$ are isolated points of $\hbox{J-}\sigma(U_p(q))$, and $$\Psi : {B}\to C(\hbox{J-}\sigma(U_p(q))\backslash\{ 0, 1\},M_2(\mathbb{C}))\oplus^{\infty} C(\hbox{J-}\sigma(p+2q)\cap\{1,2, 3\}) =: \tilde{\mathfrak{B}}$$ is a $^*$-isomorphism. Following the notation of the previous case, we have
	$$\begin{aligned}
		\tilde{p}(t) &= \Psi(p)(t) = \left(\begin{pmatrix}
			1 & 0 \\ 0 & 0
		\end{pmatrix}, (1, 0, 1)\right) \text{ and } \\ 
		\tilde{q}(t) &= \Psi(q)(t) = \left(\begin{pmatrix}
			t & \sqrt{t(1-t)} \\ \sqrt{t(1-t)} & 1-t
		\end{pmatrix}, (0, 1, 1)\right).
	\end{aligned}$$
	Defining $\tilde{w}(t) = U_{\tilde{p}}(\tilde{q})(t) = \left( \begin{pmatrix}
		t & 0 \\ 0 & 0
	\end{pmatrix}, (0, 0, 1) \right) \in \Psi(\mathfrak{B})$, we can observe that $\tilde{w}^n$ converges to $\tilde{e}_3 = \left(0, (0, 0, 1)\right)$, which (by the closed character of $\Psi(\mathfrak{B})$) also lies in this JB$^*$-subalgebra. By defining $\tilde{r}(t) = U_{\tilde{p}-\tilde{e}_3}(\tilde{q}-\tilde{e}_3)(t) = \begin{pmatrix}
		t & 0 \\ 0 & 0
	\end{pmatrix} \in \Psi_1(\mathfrak{B})$, we can then proceed in a very similar way as before in order to conclude that $\Psi(\mathfrak{B}) = \tilde{\mathfrak{B}}$.
\end{proof}

Let us note that in the proof of \Cref{p JBstar subalgebra generated by two projections} we have employed and stated some celebrated characterizations of the C$^*$-algebra generated by two projections due to G.K. Pedersen \cite[Theorem 3.4]{Ped68} and  N. Krupnik, S. Roch, and B. Silbermann \cite[Theorem 5b]{KruRothSilb}.  

A consequence of our previous proposition characterizes when two projections in an arbitrary JB$^*$-algebra are orthogonal in terms of the distance between the positive spheres of their respective Peirce $2$-subspaces. 

\begin{lemma}\label{lem:characterization-orthogonality JB*-algebra} Let $\mathfrak{A}$ be a JB$^*$-algebra. Let $p, q \in \CP (\mathfrak{A}) \setminus \{0\}$. Then, $p$ and $q$ are orthogonal if, and only if, $\|a - b\| = 1$ for all $a, b \in S_{\mathfrak{A}^+}$ such that $a = U_p(a)$, $b = U_q(b)$.
\end{lemma}

\begin{proof} Suppose first that $p \perp q$. Then, it is known that for any $a, b \in S_{\mathfrak{A}^+}$ such that $a = U_p(a)$ and $b = U_q(b)$, one has $a \perp b$. In particular, $\|a - b\| = \max\{\|a\|, \|b\|\} = 1$.
	
	In order to see the reciprocal implication, we assume that $\|a - b\| = 1$ for all $a, b \in S_{\mathfrak{A}^+}$ such that $a = U_p(a)$, and $b = U_q(b)$. Let $\mathfrak{B}$ denote the JB$^*$-subalgebra of $\mathfrak{A}$ generated by $p$ and $q$, and consider the representation $\Phi$ established in \Cref{p JBstar subalgebra generated by two projections}. Observe that, by the just quoted proposition, $p \circ q =0$ if, and only if, $\hbox{J-}\sigma(p+2q)\cap\{3\} = \emptyset$ and $\hbox{J-}\sigma(U_p(q))\cap (0,1) = \emptyset$ (the latter is equivalent to $U_p (q)$ being a projection) in cases $(a)$ and $(b)$, and $\hbox{J-}\sigma(U_p(q))\cap (0,1) = \emptyset$ in cases $(c)$ and $(d)$. We explain, for example, case $(a)$. Namely, suppose $0$ and $1$ are not accumulation points of $\hbox{J-}\sigma(U_p(q))$. Then, by \Cref{p JBstar subalgebra generated by two projections}, there exists a Jordan $^*$-isomorphism $$\Phi : \mathfrak{B}\to C(\hbox{J-}\sigma(U_p(q))\backslash\{0, 1\},S_2(\mathbb{C}))\oplus^{\infty} C(\hbox{J-}\sigma(p+2q)\cap\{1, 2, 3\})$$ satisfying \eqref{eq def of Phi1}, and hence $$\Phi (p) (t,s) =\left(\left(\begin{matrix} 1 & 0 \\
		0 & 0
	\end{matrix}\right), \Phi_2(p) (s) \right), \Phi (q) (t,s) =\left(\left(\begin{matrix} t & \sqrt{t (1-t)} \\
		\sqrt{t (1-t)} & 1-t
	\end{matrix}\right), \Phi_2(q) (s) \right),$$ for all $t\in \hbox{J-}\sigma(U_p(q))\backslash\{0, 1\}$, $s\in \hbox{J-}\sigma(p+2q)\cap\{1, 2, 3\}$, and the rest is clear.
	
	If $\hbox{J-}\sigma(p+2q)\cap\{3\} \neq \emptyset$ in case $(a)$ or $(b)$, the element $a\in \mathfrak{B}$ satisfying $\Phi (a) = (0, (0,0,1))$ in case $(a)$ and $\Phi (a) = (0, 1)$ in case $(b)$,  is a positive norm-one element. Taking $a= b \leq p,q,$ we have $\|a-b\| =0$, which contradicts our assumption.
	
	Suppose next that there exists $t_0\in \hbox{J-}\sigma(U_p(q))\cap (0,1) \neq \emptyset$ in any case of \Cref{p JBstar subalgebra generated by two projections}. Let us denote $\tilde{p} = \Phi(p)$ and $\tilde{q} = \Phi(q)$. Take $\delta > 0$ such that $\delta < \min\left\{t_0, \frac{4}{5}\right\}$, and consider the continuous function $f : [0, 1] \to \mathbb{R}$ given by
	\begin{equation*}
		f(t) = \left\{ \begin{array}{lr}
			0, & \hbox{if } 0 \leq t < \frac{\delta}{2} \\
			\frac{2}{\delta} t - 1, &  \hbox{if } \frac{\delta}{2} \leq t < \delta \\
			1, & \hbox{if } \delta \leq t < 1.
		\end{array} \right.
	\end{equation*} Clearly, the restricted function $f|_{\hbox{J-}\sigma(U_p(q))}$ {lies in}  $C(\hbox{J-}\sigma(U_p(q)))$. We can actually work with the function $f$ itself to consider the positive norm-one elements $\tilde{a} (t) = f(t) \Phi_1 (p) (t) = (f \Phi_1 (p)) (t)$ and $\tilde{b} (t) = f(t) \Phi_1 (q)  (t)$. {It} is easy to check that $\tilde{a} = U_{\tilde{p}}(\tilde{a})$ and $\tilde{b} = U_{\tilde{q}}(\tilde{b})$, and hence $a= \Phi^{-1} (\tilde{a})$ and $b= \Phi^{-1} (\tilde{b})$ are positive norm-one elements in $\mathfrak{B}\subseteq \mathfrak{A}$ with $U_p (a) = a$ and $U_q (b) = b$. Finally, by construction, we have  $$\begin{aligned}
		\|a-b\| &= \|\tilde{a}-\tilde{b}\| = \sup_{t\in \hbox{J-}\sigma(U_p(q))} f(t) \|  \Phi_1 (p)  (t) -  \Phi_1 (q)  (t)  \| \\ &= \sup_{t\in \hbox{J-}\sigma(U_p(q))} f(t) \sqrt{1-t} \leq  \sup_{t\in\left[\frac{\delta}{2},1\right]}  f(t) \sqrt{1 - t} =  \sqrt{1 - \delta} < 1,
	\end{aligned}$$ which contradicts the assumption.
\end{proof}

If in the characterization of orthogonality given in \Cref{lem:characterization-orthogonality JB*-algebra} we are interested in replacing the positive norm-one elements $a,b$ with non-zero projections we shall need to work in the setting of JBW$^*$-algebras. For this purpose we shall need a characterization of the JBW$^*$-subalgebra generated by two projections.

Let us return to the Shirshov-Cohn theorem, and its consequence affirming that the JB$^*$-subalgebra $\mathfrak{B}$ generated by two self-adjoint elements $a,b$ (and possibly $\one$) in a (possibly unital) JB$^*$-algebra $\mathfrak{A}$ is a JC$^*$-algebra (cf. \cite[Theorem 1.14]{AlfsenShultz2003}, \cite[Theorem 7.2.5]{HOS} or \cite[Proposition 2.1 and Corollary 2.2]{Wright1977}). The argument to get this conclusion is essentially the following: by a result of E.M. Alfsen, F.W. Shultz and E. St{\o}rmer \cite{AlfShulStor79GelfandNeumark} (see also \cite[Lemma 7.2.2, Theorem 7.2.3, and Corollary 4.5.8]{HOS}), if $\mathfrak{B}$ is not a JC$^*$-algebra, there exists a surjective Jordan $^*$-homomorphism $\pi: \mathfrak{B}\to H_3(\mathbb{O}^{\mathbb{C}})$ (which is automatically continuous). Since $\mathfrak{B}$ is generated, via the norm closure, by the Jordan $^*$-subalgebra algebraically generated by $a$ and $b$ (and possibly $\one$), and $H_3(\mathbb{O}^{\mathbb{C}})= \pi(\mathfrak{B})$ is clearly finite-dimensional, we conclude that $H_3(\mathbb{O}^{\mathbb{C}})= \pi(\mathfrak{B})$ must be algebraically generated by $\pi(a)$, $\pi(b)$ (and possibly $\one$) which is known to be incompatible with the fact that $H_3(\mathbb{O}^{\mathbb{C}})$ is exceptional (cf. \cite[Corollary 2.8.5 and the Shirshov-Cohn theorem]{HOS}).

Suppose now that $\mathfrak{A}$ is a JBW$^*$-algebra. What can we say about the JBW$^*$-subalgebra $\mathfrak{B}$ of $\mathfrak{A}$ generated by two self-adjoint elements $a,b$ (and possibly $\one$) in $\mathfrak{A}$? It is stated in \cite[Theorem 3.6]{vanWetering2020} that a similar argument to that given in the previous paragraph remains valid to show that $\mathfrak{B}$ is a JW$^*$-algebra, i.e. a weak$^*$-closed JB$^*$-subalgebra of some von Neumann algebra $M$. As before, if $\mathfrak{B}$ is not a JC$^*$-algebra, there exists a surjective Jordan $^*$-homomorphism $\pi: \mathfrak{B}\to H_3(\mathbb{O}^{\mathbb{C}})$. However, in this case $\mathfrak{B}$ is the weak$^*$-closure of the Jordan $^*$-subalgebra algebraically generated by $a$ and $b$ (and possibly $\one$), and we cannot guarantee the weak$^*$-continuity of the Jordan $^*$-homomorphism $\pi$. {Thus,} it is not obvious why $H_3(\mathbb{O}^{\mathbb{C}})= \pi(\mathfrak{B})$ must be algebraically generated by $\pi(a)$ {and} $\pi(b)$ (and possibly $\one$). To avoid this difficulty we include here an alternative approach. 

Let $\mathfrak{A}$ be a JBW$^*$-algebra. The $\sigma$-strong$^*$-topology (also called the strong$^*$-topology) of $\mathfrak{A}$ is the topology on $\mathfrak{A}$ defined by all the pre-Hilbert seminorms of the form $a\mapsto \varphi (a \circ a^*)^{\frac12}$ ($a\in \mathfrak{A}$), where $\varphi$ is running in the set of all positive norm-one normal functionals in $\mathfrak{A}_{*}$. The $\sigma$-strong topology of $\mathfrak{A}_{sa}$ is the topology obtained by restricting to $\mathfrak{A}_{sa}$ the $\sigma$-strong$^*$-topology of $\mathfrak{A}$, equivalently, the topology on $\mathfrak{A}_{sa}$ defined by all the pre-Hilbert seminorms of the form $a\mapsto \varphi (a^2)^{\frac12}$ ($a\in \mathfrak{A}_{sa}$), where $\varphi$ is running in the set of all positive norm-one normal functionals in $(\mathfrak{A}_{sa})_{*}$ (cf. \cite[Definition 4.1.3]{HOS} or \cite[Definition 2.3]{AlfsenShultz2003}). 

\begin{proposition}\label{p the JBW*-subalgebra gen by two projections is special} Let $\mathfrak{A}$ be a JBW$\,^*$-algebra, and let  $\mathfrak{B}$ denote the JBW$\,^*$-subalgebra of $\mathfrak{A}$ generated by two self-adjoint elements $a,b$ {\rm(}and possibly $\one${\rm)} in $\mathfrak{A}$. Then $\mathfrak{B}$ is a JW$\,^*$-algebra. 
\end{proposition}

\begin{proof} According to the arguments in the previous paragraphs, if $\mathfrak{B}$ is not a JC$^*$-algebra, there exists a surjective Jordan $^*$-homomorphism $\pi: \mathfrak{B}\to H_3(\mathbb{O}^{\mathbb{C}})$. In particular $\pi(\mathfrak{B}_{sa}) = H_3(\mathbb{O})$.
	
	As explained by E. M. Alfsen, F. W. Shultz and E. St{\o}rmer \cite[comments before Lemma 9.4]{AlfShulStor79GelfandNeumark}, although there exists no set of identities characterizing special Jordan algebras
	among all Jordan algebras, there do exist $s$-identities satisfied by all special Jordan algebras but not by all Jordan algebras. For example, as shown by C.M. Glennie (see \cite[Theorem 12, page 51]{JacobsonBook68}) by defining \begin{equation}\label{triple product} \{x,y,z\} = (x\circ y^*)\circ z + (z\circ y^*)\circ x- (x\circ z) \circ y^*,
	\end{equation}
	the $s$-identity  
	\begin{equation}\label{eq Glennie identity} \begin{aligned}
			\mathcal{G}(x,y,z) &:=2 \{x,z,x\} \circ \{y,\{z,y^2,z\},x\} - 2\{y,z,y\} \circ  \{y,\{z,y^2,z\},x\} \\
			&- \{x,\{z,\{x, \{y,z,y\} ,y\},z\},x\} + \{y, \{z, \{x, \{x,z,x\} ,y  \} ,z\}, y\}=0,
		\end{aligned}
	\end{equation} is satisfied for all $x,y,z$ in a special Jordan algebra but it is not satisfied for all $x,y,z\in  H_3(\mathbb{O}) = H_3(\mathbb{O}^{\mathbb{C}})_{sa}$. 
	
	Let $\mathfrak{B}_0$ denote the JB$^*$-subalgebra of $\mathfrak{B}$ generated by $a$, $b$ and possibly $\one$. We have already commented that $\mathfrak{B}_0$ is a JC$^*$-algebra (\cite[Theorem 7.2.5]{HOS} or \cite[Proposition 2.1 and Corollary 2.2]{Wright1977}), and hence the Glennie's $s$-identity in \Eqref{eq Glennie identity} holds for all $x,y,z\in (\mathfrak{B}_0)_{sa}$. 
	
	We claim that Glennie's $s$-identity in \Eqref{eq Glennie identity} actually holds for all $\tilde{x},\tilde{y}, \tilde{z}\in \mathfrak{B}_{sa}$. Fix $\tilde{x},\tilde{y}, \tilde{z}\in \mathfrak{B}_{sa}$.  Obviously $(\mathfrak{B}_0)_{sa}$ is weak$^*$-dense in $\mathfrak{B}_{sa}$. Observe that the weak$^*$- (called the $\sigma$-weak topology in \cite{AlfsenShultz2003} and the weak topology in \cite{HOS}) and the $\sigma$-strong topologies on $\mathfrak{B}_{sa}$ are compatible topologies (see \cite[Corollary 4.5.4]{HOS} or \cite[Proposition 2.68]{AlfsenShultz2003}), and thus the weak$^*$- and the $\sigma$-strong closures of $\mathfrak{B}_0$ in $\mathfrak{B}$ coincide. We are thus in a position to apply the Kaplansky density theorem for JBW-algebras (see \cite{AlfsenShultz2003}) to deduce the existence of three bounded nets $(x_j)_j,$ $(y_k)_k$ and $(z_i)_i$ in $(\mathfrak{B}_0)_{sa}$ converging in the $\sigma$-strong topology to $\tilde{x},\tilde{y},$ and $\tilde{z}$, respectively. It follows from the previous conclusions that $\mathcal{G}(x_j,y_k,z_i) = 0$ for all $j,k,i$. Since the Jordan product of $\mathfrak{B}$ is jointly $\sigma$-strong continuous on bounded sets (cf. \cite[Lemma 4.1.9]{HOS} or \cite[Proposition 2.4]{AlfsenShultz2003}) and the mapping $\mathcal{G}$ is given in terms of Jordan polynomials, taking $\sigma$-strong limits in $j,k,i$ in the identity $\mathcal{G}(x_j,y_k,z_i) = 0$ ($j,k,i$), we deduce that $\mathcal{G}(\tilde{x},\tilde{y}, \tilde{z}) = 0$, as claimed.
	
	Finally, since $\mathcal{G}(\tilde{x},\tilde{y}, \tilde{z}) = 0$ for all $\tilde{x},\tilde{y}, \tilde{z}\in \mathfrak{B}_{sa}$,  $\pi|_{\mathfrak{B}_{sa}} : \mathfrak{B}_{sa}\to H_3(\mathbb{O})$ is a surjective Jordan homomorphism. {Furthermore, since} $\mathcal{G}$ {can be written} in terms of Jordan products, we arrive at $\mathcal{G}(\pi(\tilde{x}),\pi(\tilde{y}), \pi(\tilde{z})) = \pi (\mathcal{G}(\tilde{x},\tilde{y}, \tilde{z})) =0,$ for all $\tilde{x},\tilde{y}, \tilde{z}\in \mathfrak{B}_{sa}.$ We have therefore shown that $H_3(\mathbb{O})$ satisfies Glennie's $s$-identity in \Eqref{eq Glennie identity}, which is impossible.
	
	The above arguments show that $\mathfrak{B}$ is a JC$^*$-algebra. Since $\mathfrak{B}$ is clearly a JBW$^*$-algebra, a standard argument in JB$^*$-algebra theory (see \cite[Corollary 2.78]{AlfsenShultz2003}) implies that $\mathfrak{B}$ is a JW$^*$-algebra. 
\end{proof}

For each self-adjoint element $a$ in a JBW$^*$-algebra $\mathfrak{A}$, there exists a smallest projection $p\in \mathfrak{A}$ satisfying $p\circ a = a$, which is called the \emph{range projection} of $a$ in $\mathfrak{A}$ and is denoted by $r(a)$ (cf. \cite[Lemma 4.2.6]{HOS}). It is known that $r(a)$ coincides with the weak$^*$-limit of a sequence of elements in the JB$^*$-subalgebra of $\mathfrak{A}$ generated by $a$. Therefore, if $\mathfrak{W}$ is any JBW$^*$-subalgebra of $\mathfrak{A}$ containing the element $a$, the range projections of $a$ in $\mathfrak{A}$ and $\mathfrak{W}$ coincide. If $\mathfrak{A}$ is a JB$^*$-algebra, we shall write $r(a)$ for the range projection of $a$ in $\mathfrak{A}^{**}$.

It is now time to study the JBW$^*$-subalgebra generated by two projections. The von Neumann subalgebra generated by two projections is a very well-known and studied object described by authors like R. Giles and H. Kummer \cite{GilesKummer71} (see also \cite[Theorem 7.1]{BottSpit2010}), and considered in celebrated monographs like \cite[Exercise 12.4.11]{KadRingVol2} and \cite[pages 306-308]{Tak}. {Next, we} state a result which gathers all the available information. Following the standard notation the symbols $\vee$ and $\wedge$ will denote the least upper bound and the greatest lower bound in the set of projections in a von Neumann algebra or in a JBW$^*$-algebra $\mathfrak{A}$ (the existence {in} the Jordan case is guaranteed by \cite[Lemma 4.2.8]{HOS}). Let $p_1,\ldots, p_k$ be a finite family in $\CP (\mathfrak{A}),$ where $\mathfrak{A}$ is a JBW$^*$-algebra. Since $\displaystyle \bigvee_{i=1}^{k} p_i = r\left(\sum_{i=1}^{k} p_i\right)$ and  $\displaystyle \bigwedge_{i=1}^{k} p_i = \one- \left(\bigvee_{i=1}^{k} (\one-p_i)\right)$,  the least upper bound and the greatest lower bound of any finite family of projections in $\mathfrak{A}$ do not change when computed in any unital JBW$^*$-subalgebra of $\mathfrak{A}$ containing this family.  

From now on, given a C$^*$-algebra $A$, we shall denote by $M_n(A)$ the C$^*$-algebra of all $n\times n$ matrices with entries in $A$ with the standard operations (cf. \cite[\S IV.3]{Tak}). The symbol $S_n(A)$ will stand for the subspace of $M_n(A)$ of all $n\times n$ symmetric matrices. If $A$ is commutative, $S_n(A)$ is actually a JC$^*$-subalgebra of $M_n(A)$.

The next proposition gathers the current knowledge on the von Neumann algebra generated by two projections.    

\begin{proposition}\label{p von Neumann algebra generated by two projecitons}$($\cite{GilesKummer71}, \cite[Theorem 7.1]{BottSpit2010}), \cite[Exercise 12.4.11]{KadRingVol2}, and \cite[pages 306-308]{Tak}$)$ Let ${A}$ be a von Neumann algebra, and let ${B}$ denote the von Neumann subalgebra of ${A}$ generated by two projections $p,q\in {A}$ and the unit element $\one$. Set $p_0 = p -p\wedge q-p\wedge q^{\perp}$ and $q_0 = q- p\wedge q- p^{\perp} \wedge q$. Let $\Gamma$ denote the set of all non-zero (mutually orthogonal) projections in $\{p\wedge q, p^{\perp}\wedge q, p\wedge q^{\perp}, p^{\perp} \wedge q^{\perp}\}$, and let $m = \sharp \Gamma\in\{0,1,2,3,4\}$. Then, for each $r\in \Gamma$ there exists a normal $^*$-homomorphism $\varphi_r : {B}\to \mathbb{C}$ satisfying $\varphi_r (r) = 1$ and $\varphi_r (e) = 0$ for all $e\in \{p\wedge q, p^{\perp}\wedge q, p\wedge q^{\perp}, p^{\perp} \wedge q^{\perp} \}\backslash\{r\}$, and the following statements hold:
	\begin{enumerate}[$(1)$]
		\item There exist a commutative von Neumann algebra $\mathcal{C}$ and positive elements $0\leq c,s\leq \one_{\mathcal{C}}$ in $\mathcal{C}$ satisfying $c^2 + s^2 =\one_{\mathcal{C}}$ (the unit in $\mathcal{C}$), the range projection of $cs$ is ${\one_{\mathcal{C}}}$, and ${\mathcal{C}}$ is generated by $c^2$ and $\one_{\mathcal{C}}$.
		\item There exists a C$\,^*$-isomorphism $\Phi: {B} \to M_2(\mathcal{C}) \oplus^{\infty} \mathbb{C}^{m},$ $\Phi (a) = (\Phi_1 (a),\Phi_2(a))$ satisfying the following properties:     
		\begin{enumerate}[$(a)$] 
			\item $\Phi_1$ vanishes on $(p_0\vee q_0)^{\perp} B (p_0\vee q_0)^{\perp}$, and the restricted mapping  $$\Phi_1|_{(p_0\vee q_0)B (p_0\vee q_0)}: (p_0\vee q_0)B (p_0\vee q_0) \to  M_2(\mathcal{C})$$ is a C$\,^*$-isomorphism with 
			$\Phi_1 (p) = \left(\begin{matrix} \one_{\mathcal{C}} & 0 \\
				0 & 0
			\end{matrix}\right),$ and $\Phi_1 (q) = \left(\begin{matrix} c^2 & c s \\
				c s & s^2
			\end{matrix}\right).$
			
			\item $\Phi_2 (a) = (\varphi_{r} (a))_{r\in \Gamma}$, and $\Phi_2 ((p_0\vee q_0)B (p_0\vee q_0)) =\{0\}$.
		\end{enumerate} 
	\end{enumerate}	We can additionally conclude that every projection in $\Gamma$ is minimal in $B$. $\hfill\Box$
\end{proposition}

Using \Cref{p von Neumann algebra generated by two projecitons}, a similar argument to that in the proof of \Cref{p JBstar subalgebra generated by two projections} gives the next description.

\begin{proposition}\label{p JBW generated by two projecitons} Let $\mathfrak{A}$ be a JBW$\,^*$-algebra, and let $\mathfrak{B}$ denote the JBW$\,^*$-subalgebra generated by two projections $p,q\in \mathfrak{A}$ and the unit element $\one$. Set $p_0 = p -p\wedge q-p\wedge q^{\perp}$ and $q_0 = q- p\wedge q- p^{\perp} \wedge q$. Let $\Gamma$ denote the set of all non-zero (mutually orthogonal) projections in $\{p\wedge q, p^{\perp}\wedge q, p\wedge q^{\perp}, p^{\perp} \wedge q^{\perp} \}$, and let $m = \sharp \Gamma\in\{0,1,2,3,4\}$. Then for each $r\in \Gamma$, there exists a normal Jordan $^*$-homomorphism $\psi_r : \mathfrak{B}\to \mathbb{C}$ satisfying $\psi_r (r) = 1$ and $\psi_r (e) = 0$ for all $e\in \{p\wedge q, p^{\perp}\wedge q, p\wedge q^{\perp}, p^{\perp} \wedge q^{\perp} \}\backslash\{r\}$, and  the following statements hold:
	\begin{enumerate}[$(1)$]
		\item There exist a commutative von Neumann algebra $\mathcal{C}$ and positive elements $0\leq c,s\leq \one_\mathcal{C}$ in $\mathcal{C}$ satisfying $c^2 + s^2 =\one_{\mathcal{C}}$ the unit in $\mathcal{C}$, the range projection of $cs$ is $\one_{\mathcal{C}}$, and ${\mathcal{C}}$ is generated by $c^2$ and $\one_{\mathcal{C}}$.
		\item There exists a Jordan $^*$-isomorphism $\Psi: \mathfrak{B} \to S_2(\mathcal{C}) \oplus^{\infty} \mathbb{C}^{m},$ $\Psi (a) = (\Psi_1 (a),\Psi_2(a))$ satisfying the following properties:     
		\begin{enumerate}[$(a)$] 
			\item $\Psi_1$ vanishes on $U_{(p_0\vee q_0)^{\perp}}(\mathfrak{B})$, and $\Psi_1|_{U_{p_0\vee q_0}(\mathfrak{B})} : U_{p_0\vee q_0}(\mathfrak{B}) \to  S_2(\mathcal{C})$ is a Jordan $^*$-isomorphism with $\Psi_1 (p) = \left(\begin{matrix} \one_{\mathcal{C}} & 0 \\
				0 & 0
			\end{matrix}\right),$ and $\Psi_1 (q) = \left(\begin{matrix} c^2 & c s \\
				c s & s^2
			\end{matrix}\right)$. 
			
			\item $\Psi_2 (a) = (\psi_{r} (a))_{r\in \Gamma}$, and $\Psi_2 (U_{p_0\vee q_0}(\mathfrak{B})) =\{0\}$.
		\end{enumerate} 
	\end{enumerate}	
\end{proposition}

\begin{proof} By \Cref{p the JBW*-subalgebra gen by two projections is special}, there exists a von Neumann algebra $A$ containing $\mathfrak{B}$ as a JW$^*$-subalgebra and both share the same unit element. Let $B$ denote the von Neumann subalgebra generated by $p,q,\one$, let $\Gamma_{_B}$ be the set of all non-zero projections in $\{p\wedge q, p^{\perp}\wedge q, p\wedge q^{\perp}, p^{\perp}\wedge q^{\perp} \}$, {and} let $m = \sharp \Gamma\in\{0,1,2,3,4\}$. Clearly, $\mathfrak{B}\subseteq B.$ Let us observe that there is no ambiguity in the elements $p\wedge q,$ $p^{\perp}\wedge q,$ $p\wedge q^{\perp},$ and $p^{\perp}\wedge q^{\perp},$ since they do not change when computed in $\mathfrak{A}$, $\mathfrak{B}$, $A,$ or $B$ (cf. the paragraph preceding \Cref{p von Neumann algebra generated by two projecitons}).  \Cref{p the JBW*-subalgebra gen by two projections is special} assures that for each $r\in \Gamma_{_B}$, there exists a normal $^*$-homomorphism $\varphi_r : {B}\to \mathbb{C}$ satisfying $\varphi_r (r) = 1$ and $\varphi_r (e) = 0$ for all $e\in \{p\wedge q, p^{\perp}\wedge q, p\wedge q^{\perp},  p^{\perp}\wedge q^{\perp} \}\backslash\{r\}$. Furthermore, by the just quoted proposition, there {exist} a commutative von Neumann algebra $\mathcal{C}$ and a C$^*$-isomorphism $\Phi: {B} \to M_2(\mathcal{C}) \oplus^\infty \mathbb{C}^{m},$ whose precise form and properties are given in the just quoted result.  Obviously, $\psi_r = \varphi_r |_{\mathfrak{B}}$ is a normal Jordan $^*$-homomorphism. 
	
	Clearly, $\mathfrak{B}\subseteq \Phi^{-1}\big(S_2(\mathcal{C}) \oplus^{\infty} \mathbb{C}^{m}\big)\subseteq B,$ where $\Phi^{-1}(S_2(\mathcal{C}) \oplus^{\infty} \mathbb{C}^{m})$	is a JBW$^*$-subalgebra of $B$ containing $p,q,$ and $\one$. We only need to show that $\mathfrak{B}= \Phi^{-1}\big(S_2(\mathcal{C}) \oplus^{\infty}  \mathbb{C}^{m}\big)$ to take $\Psi = \Phi|_{\mathfrak{B}}$. 
	
	For the sake of brevity, we shall only prove the proposition in the case $m =4$; the remaining cases are easier.
	Set $\tilde{p} = \Phi(p)$ and $\tilde{q} = \Phi(q)$. It is easy to see that $$\Phi (\mathfrak{B})\ni \tilde{b} = \tilde{q}- U_{\tilde{p}} (\tilde{q})- U_{\tilde{p}^{\perp}} (\tilde{q})  = \left(\left(\begin{matrix} 0 & cs \\
		cs & 0
	\end{matrix}\right), (0,0,0,0)\right). $$ 
	
	It can be checked that, for each natural $n$, $\Phi (\mathfrak{B})$ contains the element $$\tilde{b}^{2n-1} =  \left(\left(\begin{matrix} 0 & (cs)^{2n-1} \\
		(cs)^{2n-1} & 0
	\end{matrix}\right), (0,0,0,0)\right),$$ and hence, an application of the Stone-Weierstrass theorem as in previous arguments, leads to $\Phi (\mathfrak{B})\ni \tilde{b}^{\frac{1}{2n-1}} =  \left(\left(\begin{matrix} 0 & (cs)^{\frac{1}{2n-1}} \\
		(cs)^{\frac{1}{2n-1}} & 0
	\end{matrix}\right), (0,0,0,0)\right)$ for all $n\in \mathbb{N}$ (in terms of continuous functional calculus). It is well-known that the sequence $\left((cs)^{\frac{1}{2n-1}}\right)_n$ converges to the range projection of $cs$ in the weak$^*$-topology of $\mathcal{C}$, and thus  $\Phi (\mathfrak{B})$ contains the element $\tilde{u} = \left(\left(\begin{matrix} 0 & \one_{\mathcal{C}} \\
		\one_{\mathcal{C}} & 0
	\end{matrix}\right), (0,0,0,0)\right).$ Now, the identities 
	$$\begin{aligned}
		U_{\tilde{u}} U_{\tilde{p}} (\tilde{q}) &= 	\left(\left(\begin{matrix} 0 & 0 \\
			0 & c^2
		\end{matrix}\right), (0,0,0,0)\right) \!\in\! \Phi (\mathfrak{B}), \
		U_{\tilde{u}} U_{\tilde{p}^{\perp}} (\tilde{q}) &= 	\left(\left(\begin{matrix} s^2 & 0 \\
			0 & 0
		\end{matrix}\right), (0,0,0,0)\right) \!\in\! \Phi (\mathfrak{B}),
	\end{aligned}$$ combined with the fact that $\Phi (\mathfrak{B})$ is a JBW$^*$-algebra, and ${\mathcal{C}}$ is generated by $c^2$ and $\one_{\mathcal{C}}$ ({alternatively,} by $s^2$ and $\one_{\mathcal{C}}$), prove that $\Phi (\mathfrak{B})$ contains all the elements of the form $\left(\left(\begin{matrix} x & 0 \\
		0 & y
	\end{matrix}\right), (0,0,0,0)\right),$ with $x,y\in \mathcal{C}$. Consequently, the elements $$ \left(\left(\begin{matrix} 0 & \frac{x+y}{2} \\
		\frac{x+y}{2} & 0
	\end{matrix}\right), (0,0,0,0)\right) =  \tilde{u}\circ \left(\left(\begin{matrix} x & 0 \\
		0 & y
	\end{matrix}\right), (0,0,0,0)\right) $$ also belong to  $ \Phi (\mathfrak{B})$. The previous conclusions show that $\Phi(\mathfrak{B}) \supseteq S_2(\mathcal{C}) \oplus^\infty \{0\}$, and the rest is clear by the properties of $\Phi$.
\end{proof}

We can next establish a variant of \Cref{lem:characterization-orthogonality JB*-algebra} in terms of projections. 

\begin{lemma}\label{lem:describ-disj-proj} Let $\mathfrak{A}$ be a JBW$\,^*$-algebra, and let $p, q \in \CP (\mathfrak{A}) \setminus \{0\}$. Then, $p$ and $q$ are orthogonal if, and only if,  $\| r-s\| = 1$ for every $r,s\in \CP(\CA)\setminus \{0\}$ with $r\leq p$ and $s\leq q$.
\end{lemma}

\begin{proof} The necessity is clear since every couple of projections $r,s\in \CP(\CA)\setminus \{0\}$ with $r\leq p$ and $s\leq q$ must be orthogonal.  Let us see the ``if'' implication. Suppose that  $\|r-s\| =1$ for all $r,s\in \CP({{\CA}})\setminus \{0\}$ with $r\leq p$ and $s\leq q$. Let $\mathfrak{B}$ denote the JBW$^*$-subalgebra of $\mathfrak{A}$ generated by $p,q,$ and $\one$. Let $\mathcal{C}$ and $\Psi: \mathfrak{B} \to S_2(\mathcal{C}) \oplus^{\infty} \mathbb{C}^{m}$ be the commutative von Neumann algebra and the Jordan $^*$-isomorphism given by \Cref{p JBW generated by two projecitons}, respectively, and let $c,s\in \mathcal{C}$ be the positive elements satisfying  ${\Psi_1} (q) = \left(\begin{matrix} c^2 & c s \\
		c s & s^2
	\end{matrix}\right)$. Set $\tilde{p}=\Psi (p)$ and $\tilde{q}=\Psi (q)$. If $p\wedge q\neq 0$ we can take $0\neq r= s = p\wedge q\leq p,q$ with $\|r-s\|=0$, which is impossible. We therefore assume that $p\wedge q=0$.  Clearly the element $c$ cannot coincide with $\one_{\mathcal{C}}.$ In the case that $c=0$ (i.e., $s=\one_{\mathcal{C}}$) we have $\Psi (p)\circ \Psi (q) = 0$, equivalently, $p\perp q$ as desired. We finally deal with the case that $0\lneq s \lneq \one_{\mathcal{C}}$. Pick a non-zero projection $e\in \mathcal{C}$ such that $0 <\|e s^2\|< 1$ in $\mathcal{C}$. Consider the non-zero projection $\Psi(e) = \tilde{e} =\left(\left(\begin{matrix} e & 0 \\
		0 & e
	\end{matrix}\right), 0\right)\in S_2(\mathcal{C}) \oplus^{\infty} \mathbb{C}^{m}$. It is easy to see that $e$ operator commutes with every element in $\mathfrak{B}$ (actually $\tilde{e}$ is in the centre of the C$^*$-algebra $M_2(\mathcal{C}) \oplus^{\infty} \mathbb{C}^{m}$). Consequently $\tilde{r} := U_{\tilde{e}} (\tilde{p}) =  U_{\tilde{p}} (\tilde{e}) $ and $\tilde{s} := U_{\tilde{e}} (\tilde{q}) =  U_{\tilde{q}} (\tilde{e}) $ are two non-zero projections in $\Psi (\mathfrak{B})$ satisfying $\tilde{r}\leq \tilde{p}$ and  $\tilde{s}\leq \tilde{q}$. Finally $$\begin{aligned}
		\left\| \tilde{r}-\tilde{s} \right\|^2 &= \left\| U_{\tilde{e}} (\tilde{p}- \tilde{q}) \right\|^2 = \left\| U_{\tilde{e}} ((\tilde{p}- \tilde{q})^2) \right\|=  \left\| U_{\tilde{e}} \left(\left(\begin{matrix} s^2 & 0 \\
			0 & s^2
		\end{matrix}\right) \right) \right\| \\
		&= \left\| \left(\begin{matrix} e s^2 & 0 \\
			0 & e s^2
		\end{matrix}\right)\right\| = \| e s^2\| <1,
	\end{aligned}$$ which is also impossible.   
\end{proof}

Let us include some other consequences of our previous conclusions. 

\begin{lemma}\label{lem:describ-ord-proj} Let $\CA$ be a JBW$\,^*$-algebra, and let $p,q\in \CP (\CA) \setminus \{0\}$.
	Then $q\leq p$ if, and only if, for any $a\in \RS_{\CA^+}^{-1}$ with $\| a - p\| = 1$ we have $\| a - q\| = 1$, that is, 
	$$\sphe (p)\cap  \CA^{-1} \subseteq \sphe (q).$$
\end{lemma}

\begin{proof} Assume that $q\leq p$. Consider $a\in \RS_{\CA^+}^{-1}$ satisfying $\| a - p\| = 1$. By statements $(a)$ and $(b)$ of \Cref{lem:peak-and-zero}, there exists $\omega\in \ps(\CA)$ such that $\omega(a) = 1$ and $\omega(p) = 0$. From this, we see that $\omega(q) =0$ and hence $\omega(a-q) = 1$, which gives $\|a-q\|=1$.
	
	Conversely, suppose that the said condition holds. Let $\mathfrak{B}$ denote the JBW$^*$-subalgebra of $\CA$ generated by $p,q$ and $\one$. Consider the representation of $\mathfrak{B}$ given by \Cref{p JBW generated by two projecitons}. If $p^{\perp} \wedge q\neq 0$, we can assume $m=4$ (the other cases can be similarly treated). Since $\Psi (p) = \left(  \left(\begin{matrix} \one_{\mathcal{C}} & 0 \\
		0 & 0
	\end{matrix}\right), (1, 0, 1, 0) \right),$ and $\Psi (q) = \left(   \left(\begin{matrix} c^2 & c s \\
		c s & s^2
	\end{matrix}\right), (1, 1, 0, 0)\right)$, we can consider the element $a = \Psi^{-1} \left(    \frac12 \left(\begin{matrix} \one_{\mathcal{C}} & 0 \\
		0  & \one_{\mathcal{C}}
	\end{matrix}\right), \left(\frac12, 1, \frac12, \frac12\right) \right) \in \RS_{\CB^+}^{-1}\subseteq  \RS_{\CA^+}^{-1}$ which satisfies  $\| a- p\| = 1$ and $\|a-q \| = \frac12$, which is impossible.
	
	Assume now that $p^{\perp} \wedge q= 0$ (with $m=3$ in \Cref{p JBW generated by two projecitons}). If $q\nleq p$, $c^2 \neq \one_{\mathcal{C}}$, we can find a non-zero projection $e$ in $\mathcal{C}$ satisfying $\| e c \|<1$. Take a positive $\varepsilon$ such that $\varepsilon + \|e c\| <1$. Consider the element $$a= \Psi^{-1} \left(   \frac{\one_{\mathcal{C}}-e}{2} \left(\begin{matrix} \one_{\mathcal{C}} & 0 \\
		0  & \one_{\mathcal{C}}
	\end{matrix}\right)+ e \left(\begin{matrix} \varepsilon \one_{\mathcal{C}} & 0 \\
		0  & \one_{\mathcal{C}}
	\end{matrix}\right), \left(\frac12,\frac12, \frac12\right)\right) \in \RS_{\CB^+}^{-1}\subseteq  \RS_{\CA^+}^{-1} .$$ It is easy to check that $$1\geq  \| p-a  \| \geq \left\| \left(\begin{matrix} \one_{\mathcal{C}} & 0 \\
		0  & 0
	\end{matrix}\right) - \left(\begin{matrix} \varepsilon \one_{\mathcal{C}} & 0 \\
		0  & \one_{\mathcal{C}}
	\end{matrix}\right) \right\| \geq 1,$$ while $$\begin{aligned}
		\| q-a  \| &= \left\| \left( (\one_{\mathcal{C}}-e) \left(\begin{matrix} \frac12 \one_{\mathcal{C}}-c^2 & -cs  \\
			-cs   & \frac12 \one_{\mathcal{C}}- s^2
		\end{matrix}\right) + e \left(\begin{matrix} \varepsilon \one_{\mathcal{C}}-c^2 & -cs \\
			-cs   & \one_{\mathcal{C}}-s^2
		\end{matrix}\right), \left(\frac12, \frac12, -\frac12\right) \right) \right\| \\
		&\leq \max\left\{ \left\| \frac12 \left(\begin{matrix} \one_{\mathcal{C}} & 0 \\
			0  & \one_{\mathcal{C}} 
		\end{matrix}\right) -  \left(\begin{matrix} c^2 & cs \\
			cs  & s^2
		\end{matrix}\right) \right\| , \left\| e \left(\begin{matrix} \varepsilon \one_{\mathcal{C}} -c^2 & -cs \\
			-cs  & \one_{\mathcal{C}}-s^2
		\end{matrix}\right) \right\|, \frac12 \right\} \\
		&\leq \max\left\{\frac12,  \varepsilon \left\| \left(\begin{matrix} \one_{\mathcal{C}} & 0 \\
			0  & 0
		\end{matrix}\right) \right\| + \left\| e  \left(\begin{matrix} c^2 & c s \\
			c s & s^2- \one_{\mathcal{C}}
		\end{matrix}\right) \right\| \right\} \\
		&= \max\left\{\frac12,  \varepsilon + \left\| e \left(\begin{matrix} c^2 & c s \\
			c s & s^2-\one_{\mathcal{C}}
		\end{matrix}\right) \right\| \right\} = \max\left\{\frac12,  \varepsilon + \left\| e \left(\begin{matrix} c^2 &  c s \\
			c s & - c^2
		\end{matrix}\right) \right\| \right\} \\
		&=  \max\left\{\frac12,  \varepsilon + \left\| e c \right\| \right\}<1,
	\end{aligned}$$ where, clearly $\left\| \frac12 \left(\begin{matrix} \one_{\mathcal{C}} & 0 \\
		0  & \one_{\mathcal{C}}
	\end{matrix}\right) - \left(\begin{matrix} c^2 & c s \\
		c s & s^2
	\end{matrix}\right) \right\|=\frac12$,  and  in the penultimate equality we applied that  $$\left\| e \left(\begin{matrix} c^2 &  c s \\
		c s & - c^2
	\end{matrix}\right) \right\|^2 = \left\| e \left(\begin{matrix} c^2 &  c s \\
		c s & - c^2
	\end{matrix}\right)^2 \right\| = \left\| e \left(\begin{matrix} c^2 &  0 \\
		0 & c^2
	\end{matrix}\right) \right\| = \|e c^2\| = \| ec \|^2.$$
\end{proof}

Recall that a projection $p$ in a JBW$^*$-algebra $\CA$ is called \emph{Abelian} if $U_p (\CA)$ is an associative JBW$^*$-algebra (that is, a commutative von Neumann algebra). The representation result in \Cref{p JBW generated by two projecitons} can be also applied to characterize Abelian {projections} in terms of distances.

\begin{corollary}\label{c characterization of abelian projections}  Let $p$ be a projection in a JBW$\,^*$-algebra $\CA$. Then $p$ is Abelian if, and only if, for every $r,q\in \CP (\CA)$ with $r,q\leq p$ we have $\|q-r\|\in \{0,1\}$. 
\end{corollary}

\begin{proof} $(\Rightarrow)$ If $p$ is Abelian, every $q\in \CP (\CA)$ with $q\leq p$ is a projection in the commutative von Neumann algebra $U_p (\CA)$, and in a commutative von Neumann algebra the distance between projections is always $0$ or $1$. 
	
	$(\Leftarrow)$ Suppose now that $\|q-r\|\in \{0,1\}$ whenever $r,q\in \CP (\CA)$ with $r,q\leq p$. Since every Hermitian element in a JBW$^*$-algebra can be approximated in norm by finite linear combinations of mutually orthogonal projections, the JBW$^*$-algebra $U_p (\CA)$ will be associative if and only if any two projections in $U_p (\CA)$ operator commute. Suppose, contrary to the desired conclusion, that {$p$} is not Abelian. Then there exists $q,e\in \CP(U_p (\CA))$ which do not operator commute. If we consider the JBW$^*$-subalgebra of $U_p (\CA)$ generated by $\{q,e,p\}$ and \Cref{p JBW generated by two projecitons}, we can deduce, as in the proof of \Cref{lem:describ-ord-proj}, the existence of a projection $r\in U_p(\CA)$ with $\|q-r\|<1$, which contradicts our assumptions.  
\end{proof}

\section{Preservers of points at diametrical distance}\label{Sec: preservers of diametrical distance among projections}

As well as von Neumann algebras can be classified in terms of the Murray-von Neumann equivalence on their projection lattices, a dimension theory is also available for JBW$^*$-algebras, however the equivalence relation is somehow different and closer to unitary equivalence  (cf. \cite[\S 5]{HOS} and \cite[\S 3]{AlfsenShultz2003}). Recall that an element $s$ in a JBW$^*$-algebra $\mathfrak{A}$ is said to be a symmetry if $s=s^*$ and $s^2={\one}$. Two projections $p,q$ in $\mathfrak{A}$ are said to be \emph{equivalent} ($p \sim q$ in short) if there exist symmetries $s_1,\ldots, s_n$ in $\CA$ such that $q =  U_{s_1}U_{s_2}\cdots U_{s_n} (p)$. We also write $p \sim_n q$ in this case. If $n = 1$ we say $p$ and $q$ are \emph{exchanged by a symmetry} (cf. \cite[5.1.4]{HO1983} and \cite[Definitions 3.1 and 3.2]{AlfsenShultz2003}). 

JBW$^*$-algebras admit a classification into types $I$, $II$, and $III$ (cf. \cite[\S 5]{HOS}) in a similar fashion to that for von Neumann algebras. A projection $p\in\mathfrak{A}$ is said to be \emph{modular} if the projection lattice of $U_p (\mathfrak{A})$ is modular. The JBW$^*$-algebra $\mathfrak{A}$ is itself modular if its unit is modular.  $\mathfrak{A}$ is of type $I$ if there is an Abelian projection $p$ in $\mathfrak{A}$ whose central support, $c(p),$ is $\one$. $\mathfrak{A}$ is of type $II$ if there is a modular projection $p$ in $\mathfrak{A}$ with 
central support $\one$ and contains no non-zero Abelian projection. Finally, $\mathfrak{A}$ is of type $III$ if it contains no non-zero modular projection. Every JBW$^*$-algebra admits a unique decomposition into an orthogonal sum of JBW$^*$-algebras of type $I$, $II$ and $III$ \cite[Theorem 5.1.5]{HOS}. Each JBW$^*$-algebra of type $I$ admits a finer decomposition into summands which are either zero or JBW$^*$-algebras of type $I_n$ for $n = 0,1, 2,\ldots, \infty$ \cite[Theorem 5.3.5]{HOS}. A JBW$^*$-algebra $\mathfrak{A}$ is of \emph{type $I_n$}, where $n$ is a cardinal number, if there is a family $(p_{j})_{j\in \Gamma}$ of Abelian projections such that $c(p_j)= \one,$ $\sum_{j\in \Gamma} p_j = \one,$ and card$(\Gamma)= n$. 

Let $\mathfrak{A}$ and $\mathfrak{B}$ be JBW$^*$-algebras. A mapping $\Phi: \CP (\mathfrak{A})\to \CP (\mathfrak{B})$ is called an \emph{ortho-isomorphism} if it is a bijection preserving orthogonality in both directions. We say that $\Phi$ is an \emph{order isomorphism} if it is bijective and preserves the partial order in both directions.

\begin{proposition}\label{prop 4 point 1} Let $\CA$ and $\CB$ be JBW$\,^*$-algebras. Suppose $\Theta:\CP(\CA) \to \CP(\CB)$ is an order isomorphism
	that preserves points at diametrical distance, that is, $$\hbox{$\|p-q\| =1$ in $\CP(\CA)$ if, and only if, $\|\Theta (p) -\Theta(q) \| =1$ in $\CP(\CB)$. }$$
	Then the following statements hold: 	
	\begin{enumerate}[$(a)$] \item $\Theta$ is bi-orthogonality preserving and $\Theta(\CP({\Z(\CA)})) = \CP({\Z(\CB)})$.
		\item If $\CA$ does not contain any type $I_2$ direct summand, then $\Theta$ extends to a Jordan $^*$-isomorphism from $\CA$ onto $\CB$.
		\item Let ${e_{I_2}}\in \CP({\Z(\CA)})$ and ${f_{I_2}}\in \CP({\Z(\CB)})$ be the central projections determining the type $I_2$ parts of $\CA$ and $\CB$, respectively. Then $\Theta(e_{I_2}^{\perp}) = f_{I_2}^{\perp}$, $\Theta({e_{I_2}}) = {f_{I_2}},$ $\Theta(\U_{e_{I_2}}(\CP(\CA))) = \U_{f_{I_2}}(\CP(\CB)),$ and $\Theta(\U_{e_{I_2}^{\perp}}(\CP(\CA))) = \U_{f_{I_2}^{\perp}}(\CP(\CB))$.
	\end{enumerate}
\end{proposition}

\begin{proof}
	$(a)$ Since $\Theta$ preserves points at diametrical distance, it follows from Lemma \ref{lem:describ-disj-proj}  that $\Theta$ is bi-orthogonality preserving. Specifically, if $p$ and $q$ are orthogonal in $\CA$, then $\| r-s\| = 1$ for every $r,s\in \CP(\CA)\setminus \{0\}$ with $r\leq p$ and $s\leq q$. By the hypotheses, $\| \Theta(r)-\Theta(s)\| = 1$ for every $\Theta(r),\Theta(s)\in \CP(\CB)\setminus \{0\}$ with $\Theta(r)\leq \Theta(p)$ and $\Theta(s)\leq \Theta(q)$, and the surjectivity of $\Theta$ gives $\Theta(p) \perp \Theta(q)$ in $\CB$. Moreover, the equality $\Theta(\CP({\Z(\CA)})) = \CP({\Z(\CB)})$ follows from Lemma \ref{Lemma 2point2}$(a)$ by the same reason.

	\noindent $(b)$ Applying the reasoning from $(a)$ to $\Theta^{-1}$, it follows that $\Theta$ is an order isomorphism preserving orthogonality in both directions. Therefore $\Theta$ is an ortho-isomorphism. Theorem~4.1 in \cite{Ham21} (see also \cite{BunceWright93}) assures that $\Theta$ extends uniquely to a Jordan $^*$-isomorphism from $\CA$ onto $\CB$.
	
	\noindent $(c)$ Since $\Theta$ is an order isomorphism and $\Theta(\CP({\Z(\CA)})) = \CP({\Z(\CB)})$ (see statement $(a)$ above), it follows that $\Theta (e_{I_2}), \one-\Theta (e_{I_2}) = \Theta (e_{I_2}^{\perp})\in \CP({\Z(\CB)})$. Furthermore, $\Theta$ maps $\CP({\U_{e_{I_2}^{\perp}}(\CA)}) = \U_{e_{I_2}^{\perp}}(\CP(\CA)) = \{p\in \CP(\CA) : p\leq e_{I_2}^{\perp}\}$ onto $\CP({\U_{\Theta(e_{I_2}^{\perp})}(\CB)}) = \U_{\Theta(e_{I_2}^{\perp})}(\CP(\CB))$. Since  $\U_{e_{I_2}^{\perp}}(\CA)$ is a JBW$^*$-algebra without type $I_2$ direct summands, statement $(b)$ above assures that $\Theta|_{U_{e_{I_2}^{\perp}}(\CA)}$ extends uniquely to a Jordan $^*$-isomorphism onto $\U_{\Theta(e_{I_2}^{\perp})}(\CP(\CB))$. Therefore, $\U_{\Theta(e_{I_2}^{\perp})}(\CP(\CB))$ is a JBW$^*$-algebra without type $I_2$ direct summands, which proves that $\Theta(e_{I_2}^{\perp})\leq f_{I_2}^{\perp}$.
	By considering $\Theta^{-1}$, we arrive at $\Theta(e_{I_2}^{\perp}) = f_{I_2}^{\perp}$. Consequently, $\Theta({e_{I_2}}) = {f_{I_2}}$.
\end{proof}

It should be pointed out that the structure theory of type $I_2$ JBW$^*$-algebras differs notably from the classical classification of von Neumann algebras of type $I_2$. It is {known} that each von Neumann factor of type $I_2$ is C$^*$-isomorphic to $B(\ell_2^2)\equiv M_2 (\mathbb{C})$ (cf. \cite[Corollary V.1.28]{Tak}), and each von Neumann algebra of type $I_2$ is C$^*$-isomorphic to a von Neumann tensor product of an Abelian von Neumann algebra $\mathcal{A}$ and a von Neumann factor of type $I_2$, and therefore it can be represented in the form $C(K,M_2(\mathbb{C}))$ for some hyper-Stonean space $K$ (cf. \cite[Theorem V.1.27]{Tak}, \cite[Example 2.2.15, Theorem 2.3.3]{Sak}). It follows from the simplicity of this structure theory that two von Neumann algebras of type $I_2$ coincide if, and only if, their centres are isomorphic commutative von Neumann algebras. This is enough to guarantee that the existence of an order isomorphism preserving elements at diametrical distance between the projection lattices of two von Neumann algebras of type $I_2$ suffices to prove that they are C$^*$-isomorphic (cf. \cite[Proposition 3.2$(d)$]{LeungNgWongTAMS}).    

{If} $\CA$ is a JBW$^*$-algebra of type $I_n$, where $n$ {is} a cardinal number, then $n$ is precisely the maximal cardinality of an orthogonal family of non-zero projections in $\CA$ with the same central supports \cite[5.3.3]{HOS}. A JBW$^*$-algebra $\CA$ is a factor of type $I_2$ if, and only if, it is a spin factor (cf. \cite[Theorem 6.1.8]{HOS}, \cite[Proposition 3.37]{AlfsenShultz2003}). Here the word spin factor means complex spin factor. There exist spin factors of arbitrary {dimension} \cite[\S 6]{HOS}. To understand the representation theory of JBW$^*$-algebras of type $I_2$ we recall some definitions. 

Suppose $A$ is a commutative von Neumann algebra and $\CA$ is a JW$^*$-subalgebra of some $B(H)$. According to the usual notation, we denote by
$A \; \overline{\otimes} \; \CA$ the weak$^*$-closure of the algebraic tensor product $A\otimes \CA$ in the usual von Neumann tensor product $A
\; \overline{\otimes} \; B(H)$ of $A$ and $B(H)$. Clearly $A \; \overline{\otimes} \; \CA$ is a JBW$^*$-subalgebra of $A \; \overline{\otimes} \; B(H)$. It is well
known that every spin factor $\mathcal{V}$ is a JBW$^*$-subalgebra of some $B(H)$ (compare \cite{horn1987classification}), and thus, the von Neumann tensor product $A \; \overline{\otimes} \; \mathcal{V}$ is a JBW$^*$-algebra. Every JBW$^*$-algebra of type $I_2$ $\CA$ can be represented in the form 
$$ \CA = \bigoplus^{\ell_{\infty}}_{j\in \Gamma} A_j \; \overline{\otimes} \; \mathcal{V}_{n_j},$$ where each $A_j$ is a commutative von Neumann algebra,  each $\mathcal{V}_{n_j}$ is a spin factor with dim$(\mathcal{V}_{n_j}) ={n_j}$ and the $n_j$'s are mutually different cardinal numbers (see \cite[Classification Theorem 1.7, $(1.10)$, and Theorem 4.1]{horn1987classification} and \cite[Theorem 2 and subsequent comments]{Stacey82}). Clearly, $\CA$ is uniquely determined by the $A_j$'s and the cardinal numbers $n_j$'s. For example, the type $I_2$ von Neumann factor $M_2(\mathbb{C})$, regarded as a JBW$^*$-algebra with its Jordan product, corresponds to the four-dimensional spin factor $\mathcal{V}_4$.

It remains to understand the notion of spin factor. One of the equivalent definitions of spin factor reads as follows: A (complex) \emph{spin factor} is a complex Hilbert space $\mathcal{V}$ (with dimension $\geq 3$) together with a conjugation (i.e., a conjugate-linear isometry of {period $2$}) $x\mapsto \overline{x}$, and a distinguished norm-one element $\one = \overline{\one}\in \mathcal{V}$ with the Jordan product $$a\circ b= a\circ_{\one} b= \{a,\one, b\} = \langle a | \one \rangle b + \langle b | \one \rangle a - \langle a | \overline{b} \rangle \one \ \ (a,b\in \mathcal{V}),$$ involution $$a^* = \{\one,a,\one\}= 2 \langle \one | a \rangle \one - \overline{a} \ \ (a\in \mathcal{V}),$$ and the norm $$\| a \|^2 =  \langle a | a \rangle + \left( \langle a | a \rangle^2 - | \langle a | \overline{a}\rangle |^2 \right)^{\frac12},$$ where the triple product $\{\cdot,\cdot,\cdot\}$ is the one we considered in \Eqref{triple product}. The interested reader can consult the references \cite{Harris1974,Harris1981,HervIsi1992}, \cite[\S 3]{FriedmanBook}, and \cite[\S 1.4]{EsPeVill24} for a more detailed introduction and connections with other reformulations. 

A \emph{real spin factor} is, by definition, a JBW-algebra which coincides with the self-adjoint part $\mathcal{V}_{sa} = \{a\in \mathcal{V} : a^{*} = a \}$ of a complex spin factor $\mathcal{V}$ (see \cite[\S 6]{HO1983}). The space $\mathcal{V}^{-} :=\{x\in \mathcal{V} : \overline{x} =x\}$ is a real Hilbert space whose inner product is given by $\Re\hbox{e} \langle x | y \rangle = \langle x | y \rangle$ ($x,y\in \mathcal{V}^{-}$). It is well-known that spin factors have {rank $2$}, i.e., the maximum number of mutually orthogonal non-zero projections in a spin factor is $2$ (they actually have {rank $2$} as JB$^*$-triples \cite[page 210]{Ka97}). It follows that every non-trivial projection must be minimal, and in particular a minimal tripotent (i.e. an element $e$ satisfying $\{e,e,e\}= e$). A non-zero tripotent $e$ is called \emph{minimal} if its Peirce-$2$ subspace $\mathcal{V}_2(e) = \{x\in \mathcal{V} : \{e,e,x\} = x\}$ coincides with $\mathbb{C} e$. It is known that all minimal tripotents in $\mathcal{V}$ are of the form $\frac{a+\mathrm{i} b}{2}$, where $a,b$ are two norm-one elements in $\mathcal{V}^{-}$ with $\langle a | b \rangle =0$ (see, for example, \cite[Lemma 6.1]{HamKalPe2020}). It is a straight consequence of the above that the set of projections in $\mathcal{V}$ reduces to $\{0,\one\}$ and the minimal projections of the form $\frac{\one+\mathrm{i} b}{2}$ with $b$ in the unit sphere of $\mathcal{V}^{-}$ and $\langle \one | b\rangle =0$ (this decomposition is unique). In other words, \begin{equation}\label{eq projections in spin} \CP (\mathcal{V}) = \{0, \one \}\cup \left\{ \frac{\one+\mathrm{i} b}{2} : b\in \mathcal{V}^{-}, \ \|b\| =1, \ \langle \one | b\rangle =0 \right\}.
\end{equation} We write $\CP_{\min} (\mathcal{V})$ for the minimal projections in $\mathcal{V}$. Let us write $\{\one\}^{\perp}_{\mathcal{V}^{-}}$ for the orthogonal complement of $\one$ in the real Hilbert space $\mathcal{V}^{-}$, that is  $\{\one\}^{\perp}_{\mathcal{V}^{-}} = \{b\in \mathcal{V}^{-} : \langle \one | b\rangle =0 \}$. We have established a bijection \begin{equation}\label{eq bijection Upsilon} \Upsilon_{_\mathcal{V}}: \CP_{\min} (\mathcal{V}) \longrightarrow \RS_{\{\one\}^{\perp}_{\mathcal{V}^{-}}}, \ \ p=\frac{\one+\mathrm{i} b}{2} \mapsto \Upsilon_{_\mathcal{V}} (p)=  b.
\end{equation} Furthermore, the identity \begin{equation}\label{eq distance between minimal projections spin} \|p-q\| =  \left\| \frac{ \one + \mathrm{i}\Upsilon_{_\mathcal{V}} (p)}{2} - \frac{ \one + \mathrm{i}\Upsilon_{_\mathcal{V}} (q)}{2} \right\|= \frac12 \| \Upsilon_{_\mathcal{V}} (p) - \Upsilon_{_\mathcal{V}} (q)\|,
\end{equation} shows that $\frac12 \Upsilon_{_\mathcal{V}}$ is an isometry onto $\RS_{\{\one\}^{\perp}_{\mathcal{V}^{-}}}$. It is easy to see from \Eqref{eq distance between minimal projections spin} that for every $p,q\in \CP_{\min} (\mathcal{V})$ we have $\|p-q\|=1,$ if and only if, $\Upsilon_{_\mathcal{V}} (p) = -\Upsilon_{_\mathcal{V}} (q)$ if, and only if, $q = \one-p$.

It is an open problem whether, like in the case of type $I_2$ von Neumann algebras, the existence of an order isomorphism preserving orthogonality in both directions between the lattices of projections of two JBW$^*$-algebras of type $I_2$ implies that these two JBW$^*$-algebras are Jordan $^*$-isomorphic. The next example shows that we cannot expect a positive answer to this question.

\begin{example}\label{counterexample preservers of diametrical distance spin} Let us consider the spin factors $\mathcal{V}_3$ and $\mathcal{V}_4$, which are not Jordan $^*$-isomorphic. In this case $\{\one_{_{\mathcal{V}_3}}\}^{\perp}_{\mathcal{V}_3^{-}}\equiv \RS_{\ell_2^2(\mathbb{R})}$ and $\{\one_{_{\mathcal{V}_4}}\}^{\perp}_{\mathcal{V}_4^{-}}\equiv \RS_{\ell_2^3(\mathbb{R})}$. Take two subsets $\mathcal{S}_2\subset \RS_{\ell_2^2(\mathbb{R})}$ and $\mathcal{S}_3\subset \RS_{\ell_2^3(\mathbb{R})}$ with the properties $\mathcal{S}_k\cap \left(-\mathcal{S}_k\right)= \emptyset$ and $\mathcal{S}_k\cup \left(-\mathcal{S}_k\right) = \RS_{\ell_2^k(\mathbb{R})}$ for $k=2,3$. Since $\mathcal{S}_2$ and $\mathcal{S}_3$ have the same cardinality, we can find a bijection $g:\mathcal{S}_2\to \mathcal{S}_3$, which is extended to a bijection $\tilde{g}:  \RS_{\ell_2^2(\mathbb{R})}\to \RS_{\ell_2^3(\mathbb{R})}$ by $\tilde{g} (-x) = -g(x)$ for all $x\in \mathcal{S}_2$. Composing with $\Upsilon_{_{\mathcal{V}_3}}$ and $\Upsilon_{_{\mathcal{V}_4}}$ we can find a bijection $\Theta: \CP_{\min} (\mathcal{V}_3) \to \CP_{\min} (\mathcal{V}_4)$, which is extended to $\{0,\one_{_{\mathcal{V}_3}}\}$ by $\Theta (0)=0$ and $\Theta (\one_{_{\mathcal{V}_3}}) = \one_{_{\mathcal{V}_4}}$. We have defined a bijection $\Theta: \CP (\mathcal{V}_3) \to \CP (\mathcal{V}_4).$
	
	Given $p,q\in  \CP (\mathcal{V}_3)$ we have $p\leq q$ if, and only if, $p=0$, $q\in \CP_{\min} (\mathcal{V}_3)$ or $p\in \CP_{\min} (\mathcal{V}_3)$ and $q = \one_{_{\mathcal{V}_3}}$, it is easy to see that $\Theta (p)\leq \Theta (q)$ in all the cases. This shows that $\Theta$ is an order isomorphism. 
	
	We have also seen after the formula for the distance in \Eqref{eq distance between minimal projections spin}, $\|p-q\|=1$ for some $p,q\in \CP (\mathcal{V}_3)$ if, and only if, one of the next cases holds: $p= 0$, $q\in \CP_{\min} (\mathcal{V}_3)\cup \{\one_{_{\mathcal{V}_3}}\},$ or $p\in \CP_{\min} (\mathcal{V}_3),$ $q= \one_{_{\mathcal{V}_3}}$, or $p,q\in \CP_{\min} (\mathcal{V}_3)$ with $q = \one_{_{\mathcal{V}_3}}-p$. In the first two cases, we clearly have $\|\Theta (p) -\Theta (q)\| = 1$. In the remaining case, we observe that, by construction, $$\Theta (q) = \frac{\one_{_{\mathcal{V}_4}}+ \mathrm{i}\Upsilon_{_{\mathcal{V}_4}} (q)}{2}= \frac{\one_{_{\mathcal{V}_4}}+ \mathrm{i}\Upsilon_{_{\mathcal{V}_4}} (\one_{_{\mathcal{V}_3}}-p)}{2} = \frac{\one_{_{\mathcal{V}_4}}- \mathrm{i}\Upsilon_{_{\mathcal{V}_4}} (p)}{2} = \one_{_{\mathcal{V}_4}}- \Theta (p),$$ and thus $\|\Theta (p) -\Theta (q)\| = 1.$ Therefore $\Theta$ preserves points at diametrical distance.   
\end{example}

Despite the difficulties arising from the previous counterexample, if we replace ``preservation of points at diametrical distance in both directions'' with  ``preservation of points at distance $\frac{\sqrt{2}}{2}$ in both directions'', we can get a positive answer to the problem posed above in the case of atomic JB$^*$-algebras of type $I_2$ (i.e., $\ell_{\infty}$-sums of spin factors, i.e., $A_j$ is atomic for all $j\in \Gamma$ in the standard decomposition).

Mappings preserving points at distance $1$ in both directions are usually said to have the \emph{strong distance one preserving property} (SDOPP, for short) as can be seen in the survey paper \cite{Ding2006}. According to this terminology, we shall say that a mapping preserving points at distance $\frac{\sqrt{2}}{2}$ in both directions has the \emph{strong distance-$\frac{\sqrt{2}}{2}$ preserving property}.

\begin{proposition}\label{p atomic type I2} Let $\CA$ and $\CB$ be two atomic JBW$\,^*$-algebras of type $I_2$. Suppose there exists an order isomorphism $\Theta: \CP (\CA)\to \CP(\CB)$ satisfying the strong distance-$\frac{\sqrt{2}}{2}$ preserving property, that is, $$\|p-q\| = \frac{\sqrt{2}}{2} \Leftrightarrow \|\Theta (p)-\Theta (q) \| = \frac{\sqrt{2}}{2}.$$ Then, $\Theta$ preserves orthogonality in both directions, and the JBW$\,^*$-algebras $\CA$ and $\CB$ are Jordan $^*$-isomorphic. Furthermore, if $\Theta$ is actually an isometry, we can find a unique extension of $\Theta$ to a Jordan $^*$-isomorphism from $\CA$ onto $\CB$. 
\end{proposition}

\begin{proof} Let $\Theta : \CP (\CA)\to \CP(\CB)$ be an order isomorphism preserving points at distance $\frac{\sqrt{2}}{2}$. Let us fix some notation. By the hypotheses there exist two families $\{\mathcal{V}_{n_i}\}_{i\in \Gamma_1}$ and $\{\mathcal{W}_{n_j}\}_{j\in {\Gamma_2}},$ where  $\mathcal{V}_{n_i}$ and $\mathcal{W}_{n_j}$ are complex spin factors with dimensions $n_i$ and $n_j$, respectively, $\displaystyle \CA= \bigoplus_{i\in \Gamma_1}^{\ell_{\infty}} \mathcal{V}_{n_i}$ and $\displaystyle \CB= \bigoplus_{j\in {\Gamma_2}}^{\ell_{\infty}} \mathcal{W}_{n_j}$. We write $\one_i$ and $\one_j$ for the unit in $\mathcal{V}_{n_i}$ and $\mathcal{W}_{n_j}$, respectively. Note that a pair of non-trivial projections $p,q\in \mathcal{V}_{n_i}$ are orthogonal if, and only if, they both are minimal (in $\CA$ or in $\mathcal{V}_{n_i}$) and $p+q = \one_i$. 
	
	Let us take a non-zero projection $p\lneq  \one_i$. Clearly, $p$ must be minimal in $ \mathcal{V}_{n_i}$ and in $\CA$, and hence of the form $p =\frac{\one_i+ \rm{i} b}{2}\in \CA$ with $\overline{b}=b\in \mathcal{V}_{n_i}$, $\|b\| =1$, $\langle \one_i | b\rangle =0$. Projections in $\CA$ are elements of the form $p= (p_i)_{i\in \Gamma_1}$, where $p_i$ is a projection in $\mathcal{V}_{n_i}$ for all $i\in \Gamma_1$. Note that given $p,q\in \mathcal{P} (\CA)$ with $p_j \neq 0 = q_j$ for some $j\in \Gamma_1,$ we have $\|p-q\|=1>\frac{\sqrt{2}}{2}$. It is not hard to check from this observation, together with \eqref{eq bijection Upsilon} and \eqref{eq distance between minimal projections spin}, that 
	\begin{align*}
		\operatorname{Sph}_{_{\mathcal{P}(\CA)}}^{\frac{\sqrt{2}}{2}} \left(p\right) & := \left\{ q\in \mathcal{P} (\CA) : \|p-q\| = \frac{\sqrt{2}}{2} \right\} = \left\{ q\in \mathcal{P} (\mathcal{V}_{n_i}) : \|p-q\| = \frac{\sqrt{2}}{2} \right\} \\
		&= \left\{ \frac{\one_i+ \mathrm{i} c}{2} \in \mathcal{P} (\mathcal{V}_{n_i}) : \overline{c}=c, \|c\| =1, \langle \one_i | c\rangle= \langle b | c\rangle =0  \right\} \\ & = \left\{ q\in \mathcal{P} (\mathcal{V}_{n_i}) : \|(\one_i - p)-q\| = \frac{\sqrt{2}}{2} \right\} \\
		&=  \left\{ q\in \mathcal{P} (\CA) : \|(\one_i - p)-q\| = \frac{\sqrt{2}}{2} \right\} = \operatorname{Sph}_{_{\mathcal{P}(\CA)}}^{\frac{\sqrt{2}}{2}} \left(\one_i - p\right),
	\end{align*} and, since for each norm-one element $\xi$ in a real Hilbert space $H$ the Euclidean bi-orthogonal complement of $\{\xi\}$ is $\mathbb{R} \xi$, we have \begin{equation}\label{eq double sphere minimal under onei}  \operatorname{Sph}_{_{\mathcal{P}(\CA)}}^{\frac{\sqrt{2}}{2}} \left( \operatorname{Sph}_{_{\mathcal{P}(\CA)}}^{\frac{\sqrt{2}}{2}} (p)  \right)  = \left\{ r \in \mathcal{P}(\CA) : \|r - q\| = \frac{\sqrt{2}}{2} \text{ for all } q \in \operatorname{Sph}_{_{\mathcal{P}(\CA)}}^{\frac{\sqrt{2}}{2}} (p)   \right\}= \Big\{ p, \one_i -p \Big\}.
	\end{equation} Clearly, \begin{equation}\label{eq double sphere around onei or zero} \operatorname{Sph}_{_{\mathcal{P}(\CA)}}^{\frac{\sqrt{2}}{2}} \left(\one_i\right)=\operatorname{Sph}_{_{\mathcal{P}(\CA)}}^{\frac{\sqrt{2}}{2}} \left(0\right)=\emptyset.
	\end{equation}
	
	Take now $p_1\in \CP_{\min}(\mathcal{V}_{n_{i_1}})$ and $p_2\in \CP_{\min}(\mathcal{V}_{n_{i_2}})$, with $i_1\neq i_2$. We claim that
	\begin{equation}\label{eq:double-sphere-two-atoms-type-I2}
		\operatorname{Sph}_{_{\CP(\CA)}}^{\frac{\sqrt{2}}{2}}
		\left(
		\operatorname{Sph}_{_{\CP(\CA)}}^{\frac{\sqrt{2}}{2}}(p_1+p_2)
		\right)
		=\{p_1+p_2\}.
	\end{equation}
	To see this, we begin by observing that the arguments employed above also show that $$\operatorname{Sph}_{_{\CP(\CA)}}^{\frac{\sqrt{2}}{2}}(p_1+p_2) = \left\{ q_1+q_2: q_k\in \CP_{\min}(\mathcal{V}_{n_{i_k}}), \ \max\Big\{ \|q_1-p_1\|, \|q_2-p_2\|\Big\} = \frac{\sqrt{2}}{2}  \right\}.$$ Note that elements in $\operatorname{Sph}_{_{\CP(\CA)}}^{\frac{\sqrt{2}}{2}}(p_1+p_2)$ must have zero components in all factors $\mathcal{V}_{n_{i}}$ with $i\neq i_1,i_2$, since otherwise the distance to $p_1+p_2$ would be $1$. We claim now that the identity \begin{equation}\label{eq double sphere radious1sqrt2 p1+p2} \operatorname{Sph}_{_{\CP(\CA)}}^{\frac{\sqrt{2}}{2}}\left(\operatorname{Sph}_{_{\CP(\CA)}}^{\frac{\sqrt{2}}{2}}(p_1+p_2)\right) := \left\{ r\in \mathcal{P} (\CA) : \|r-q\| = \frac{\sqrt{2}}{2},\ \forall q\in \operatorname{Sph}_{_{\CP(\CA)}}^{\frac{\sqrt{2}}{2}}(p_1+p_2)\right\} = \left\{p_1+p_2\right\}
	\end{equation} holds. Namely, let $r$ be a projection in $\CA$ at distance $\frac{\sqrt{2}}{2}$ from all the elements in $\operatorname{Sph}_{_{\CP(\CA)}}^{\frac{\sqrt{2}}{2}}(p_1+p_2)$. As above, all components of $r$ outside $\mathcal{V}_{n_{i_1}}\oplus \mathcal{V}_{n_{i_2}}$ must be zero, and its components in these two factors must be minimal projections. Having in mind that $\| r -q\|= \frac{\sqrt{2}}{2}$ for all $q = p_1+q_2$, with $\|q_2-p_2\|=\frac{\sqrt{2}}{2}$ and all $q=q_1+p_2$, with $\|q_1-p_1\|=\frac{\sqrt{2}}{2}$, we deduce that $r = p_1 +p_2$ as desired. The formulae above also hold in $\CB$.
	
	We shall next show that for each $i\in \Gamma_1$, there exists a unique $\sigma (i)\in \Gamma_2$ satisfying $\Theta (\one_i) = \one_{\sigma(i)}$. Namely, $\Theta(\one_i)$ is clearly a non-zero projection in $\displaystyle \CB= \bigoplus_{j\in {\Gamma_2}}^{\ell_{\infty}} \mathcal{W}_{n_j}$. If $\Theta(\one_i)$ admits two non-zero components in two different factors $\mathcal{W}_{n_{j_1}}$ and $\mathcal{W}_{n_{j_2}}$ with $j_1\neq j_2$, we can find minimal projections $q_1\in \mathcal{W}_{n_{j_1}}$ and $q_2 \in \mathcal{W}_{n_{j_2}}$ such that $q_1 +q_2\leq \Theta (\one_i)$. It follows from \eqref{eq double sphere radious1sqrt2 p1+p2} that $\operatorname{Sph}_{_{\CP(\CB)}}^{\frac{\sqrt{2}}{2}}\left(\operatorname{Sph}_{_{\CP(\CB)}}^{\frac{\sqrt{2}}{2}}(q_1+q_2)\right) =\{q_1+q_2\}$. 
	
	By the hypotheses on $\Theta$, there exists a non-zero projection $p_3\leq \one_i$ with $\Theta (p_3) = q_1+q_2$. If $p_3 = \one_i$ we have $\operatorname{Sph}_{_{\mathcal{P}(\CA)}}^{\frac{\sqrt{2}}{2}} \left(p_3\right)=\emptyset$ (see \eqref{eq double sphere around onei or zero}) and $\operatorname{Sph}_{_{\mathcal{P}(\CA)}}^{\frac{\sqrt{2}}{2}} \left( \operatorname{Sph}_{_{\mathcal{P}(\CA)}}^{\frac{\sqrt{2}}{2}} \left(p_3\right) \right)= \operatorname{Sph}_{_{\mathcal{P}(\CA)}}^{\frac{\sqrt{2}}{2}} \left(\emptyset\right) = \mathcal{P}(\CA)$,  while if $p_3\neq \one_i$ we have $\operatorname{Sph}_{_{\mathcal{P}(\CA)}}^{\frac{\sqrt{2}}{2}} \left( \operatorname{Sph}_{_{\mathcal{P}(\CA)}}^{\frac{\sqrt{2}}{2}} (p_3)  \right) = \Big\{ p_3, \one_i -p_3 \Big\}$ (cf. \eqref{eq double sphere minimal under onei}). We therefore have $$\{ \sharp \mathcal{P}(\CA) ,2\}\ni \sharp \Theta \left( \operatorname{Sph}_{_{\mathcal{P}(\CA)}}^{\frac{\sqrt{2}}{2}} \left( \operatorname{Sph}_{_{\mathcal{P}(\CA)}}^{\frac{\sqrt{2}}{2}} (p_3)  \right) \right) = \sharp \left( \operatorname{Sph}_{_{\mathcal{P}(\CB)}}^{\frac{\sqrt{2}}{2}} \left( \operatorname{Sph}_{_{\mathcal{P}(\CB)}}^{\frac{\sqrt{2}}{2}} (q_1+q_2)  \right) \right) =1, $$ which is impossible. So, there exists a unique $\sigma (i)\in \Gamma_2$ such that $\Theta (\one_i) \in \mathcal{W}_{n_{\sigma(i)}}$. By comparing \eqref{eq double sphere minimal under onei} and \eqref{eq double sphere around onei or zero} we arrive at $\Theta (\one_i) = \one_{\sigma(i)}$. Applying the same argument to $\Theta^{-1}$, we see that every summand unit $\one_j$ of $\CB$ is of the form $\Theta(\one_i)$. Hence $\sigma:\Gamma_1\to\Gamma_2$ is a bijection.
	
	It now follows from the hypotheses on $\Theta$ that 
	\begin{equation}\label{eq theta maps factors to factors} \begin{aligned}
			\CP (\mathcal{W}_{n_{\sigma(i)}}) &= \left\{ r\in \CP(\CB) : r\leq \Theta(\one_i) = \one_{\sigma(i)} \right\} = \Theta \Big( \big\{ e\in \CP(\CA) : e\leq \one_i \big\} \Big) = \Theta \left( \CP (\mathcal{V}_{n_i}) \right).
		\end{aligned}
	\end{equation}
	
	We claim that $\Theta$ preserves orthogonality in both directions. Having in mind \eqref{eq theta maps factors to factors}, we can reduce to the case in which $\CA$ and $\CB$ are spin factors, where a pair of non-trivial projections $(p,q)$ are orthogonal if and only if they both are minimal and $p+q = \mathbf{1}$. Let us take a minimal projection $p =\frac{\mathbf{1}+ \rm{i} b}{2}\in \CA$ with $\overline{b}=b$, $\|b\| =1$, $\langle \mathbf{1} | b\rangle =0$. It follows from \eqref{eq double sphere minimal under onei} and the hypotheses on $\Theta$ that $$\Theta(\mathbf{1} - p) \in \operatorname{Sph}_{_{\mathcal{P}(\CB)}}^{\frac{\sqrt{2}}{2}} \left( \operatorname{Sph}_{_{\mathcal{P}(\CB)}}^{\frac{\sqrt{2}}{2}} (\Theta(p))  \right) = \Big\{ \Theta(p), \mathbf{1} - \Theta(p) \Big\},$$ and thus $\Theta(\mathbf{1} - p) = \mathbf{1} - \Theta(p) \perp \Theta(p)$. The same argument applied to $\Theta^{-1}$ ensures that $\Theta$ preserves orthogonality inside each factor in both directions. Note that, as a consequence, $\Theta$ preserves pairs of projections at distance $1$ on each factor.
	
	Moreover, let $p=(p_i)_{i\in\Gamma_1},q=(q_i)_{i\in\Gamma_1}\in\CP(\CA)$ with $p\perp q$. Since $\CA=\bigoplus_{i\in\Gamma_1}^{\ell_\infty}\mathcal V_{n_i},$ orthogonality is determined componentwise. Thus, for each $i\in\Gamma_1$, either $p_i=0$, or $q_i=0$, or $p_i,q_i$ are two minimal projections in $\mathcal V_{n_i}$ satisfying $p_i+q_i=\one_i$. Since $\Theta(\one_i)=\one_{\sigma(i)}$ for all $i \in \Gamma_1$, if $p_i$ is minimal, then $\Theta(\one_i-p_i)=\one_{\sigma(i)}-\Theta(p_i)$. Therefore, on each component, $\Theta(p)_{\sigma(i)}\perp \Theta(q)_{\sigma(i)}$. Since projections belonging to different summands of $\CB$ are automatically orthogonal, it follows that $\Theta(p)\perp\Theta(q)$. Applying the same argument to $\Theta^{-1}$, we conclude that $\Theta$ preserves orthogonality in both directions.
	
	Let us now prove that $\CA$ and $\CB$ are Jordan $^*$-isomorphic. Fix $i\in \Gamma_1$, and take any two norm-one elements $b,c\in \{\one_i\}^{\perp}_{\mathcal{V}_{n_i}^{-}}$ with $\langle b|c \rangle =0$. The projections $p_b=\frac{\one_i+\mathrm{i} b}{2}$ and $p_c = \frac{\one_i +\mathrm{i} c}{2}$ are at distance $\frac{\sqrt{2}}{2}$, and thus $\left\|\Theta(p_b) -\Theta(p_c) \right\|=\frac{\sqrt{2}}{2}$, which proves that $$\left\| \Upsilon_{_{\mathcal{W}_{n_{\sigma(i)}}}} (\Theta(p_b)) - \Upsilon_{_{\mathcal{W}_{n_{\sigma(i)}}}} (\Theta(p_c))\right\| = {\sqrt{2}}.$$ This condition implies that $$\left\langle \Upsilon_{_{\mathcal{W}_{n_{\sigma(i)}}}} (\Theta(p_b)) \Big| \Upsilon_{_{\mathcal{W}_{n_{\sigma(i)}}}} (\Theta(p_c))\right\rangle =0, \hbox{ in } \{\one_{\sigma(i)}\}^{\perp}_{\mathcal{W}_{n_{\sigma(i)}}^{-}}.$$ We have therefore shown that $n_i= \text{dim}(\mathcal{V}_{n_i})\leq \text{dim}(\mathcal{W}_{n_{\sigma(i)}}) = n_{\sigma(i)}$.  Applying the same argument to $\Theta^{-1}$ we derive that $n_i= n_{\sigma(i)}$. The final conclusion follows from the arbitrariness of $i\in \Gamma_1$. 
	
	To see the last statement, suppose that $\Theta$ is an isometry. As seen before, for each $i\in \Gamma_1$ there exists a unique $\sigma(i)\in \Gamma_2$ such that $\Theta(\one_{_{\mathcal{V}_{n_i}}}) = \one_{_{\mathcal{W}_{n_{\sigma(i)}}}},$ $n_i = n_{\sigma(i)}$ as cardinal numbers, and a bijection $$T_i=  \Upsilon_{\mathcal{W}_{n_{\sigma(i)}}}\circ  \Theta|_{\CP_{\min} (\mathcal{V}_{n_i}) }\circ  \Upsilon_{\mathcal{V}_{n_i}}^{-1} : \RS_{\{\one_i\}^{\perp}_{\mathcal{V}_{n_i}^{-}}} \longrightarrow \RS_{\{\one_{\sigma(i)}\}^{\perp}_{\mathcal{W}_{n_{\sigma(i)}}^{-}}}.$$ By the new hypotheses on $\Theta$ and \Eqref{eq distance between minimal projections spin} we have $$\begin{aligned}
		\|T_i (b) -T_i (c)\| &= \left\| \Upsilon_{\mathcal{W}_{n_{\sigma(i)}}}\circ  \Theta\circ  \Upsilon_{\mathcal{V}_{n_i}}^{-1} (b) -\Upsilon_{\mathcal{W}_{n_{\sigma(i)}}}\circ  \Theta\circ  \Upsilon_{\mathcal{V}_{n_i}}^{-1} (c) \right\| \\
		&= 2 \left\| \Theta\circ  \Upsilon_{\mathcal{V}_{n_i}}^{-1} (b) - \Theta\circ  \Upsilon_{\mathcal{V}_{n_i}}^{-1} (c) \right\| = 2 \left\|  \Upsilon_{\mathcal{V}_{n_i}}^{-1} (b) - \Upsilon_{\mathcal{V}_{n_i}}^{-1} (c) \right\| = \|b-c\|,
	\end{aligned}
	$$ for all $b,c\in \RS_{\{\one_i\}^{\perp}_{\mathcal{V}_{n_i}^{-}}}$, witnessing that $T_i$ is a surjective isometry between the unit spheres of the real Hilbert spaces $\{\one_i\}^{\perp}_{\mathcal{V}_{n_i}^{-}}$ and $\{\one_{\sigma(i)}\}^{\perp}_{\mathcal{W}_{n_{\sigma(i)}}^{-}}$. By the solution to Tingley's problem for Hilbert spaces provided by G.G. Ding \cite[Corollary 2]{Ding2002}, there {exists} an extension of $T_i$ to a surjective real linear isometry $\widehat{T}_i: \{\one_i\}^{\perp}_{\mathcal{V}_{n_i}^{-}}\longrightarrow \{\one_{\sigma(i)}\}^{\perp}_{\mathcal{W}_{n_{\sigma(i)}}^{-}},$ which can be straightforwardly extended to a surjective real linear isometry denoted by the same symbol $\widehat{T}_i: \mathcal{V}_{n_i}^{-}\longrightarrow \mathcal{W}_{n_{\sigma(i)}}^{-}$ with the additional property that $ \widehat{T}_i (\one_i) = \Theta({\one_i}) = \one_{\sigma(i)}$. It is known (cf. \cite{HervIsi1992}), and easy to check, that the mapping $\Phi_i:  \mathcal{V}_{n_i} = \mathcal{V}_{n_i}^{-}\oplus \mathrm{i} \mathcal{V}_{n_i}^{-} \to \mathcal{W}_{n_{\sigma(i)}}= \mathcal{W}_{n_{\sigma(i)}}^{-}\oplus \mathrm{i} \mathcal{W}_{n_{\sigma(i)}}^{-}$, $\Phi_i (x+\mathrm{i} y ) = \widehat{T}_i (x) + \mathrm{i} \widehat{T}_i (y)$ is a triple isomorphism for the triple product given in \Eqref{triple product}. By construction, $\Phi_i (\one_i) = \Theta({\one_i})= \one_{\sigma(i)}$. Since the Jordan product and involution are given in terms of triple products and the unit element (i.e., $x\circ y = \{x,\one, y\}$ and $x^* = \{\one, x, \one\}$), the mapping $\Phi_i$ is a Jordan $^*$-isomorphism. By definition, for each minimal projection $p = \frac{\one_i+\mathrm{i} \Upsilon_{\mathcal{V}_{n_i}} (p)}{2}$ in $\mathcal{V}_{n_i}$ we have $$ \Phi_i (p) = \frac{\widehat{T}_i(\one_i)+\mathrm{i} \widehat{T}_i(\Upsilon_{\mathcal{V}_{n_i}} (p))}{2} = \frac{\one_{\sigma(i)}+\mathrm{i} \Upsilon_{\mathcal{W}_{n_{\sigma(i)}}} (\Theta (p))}{2} = \Theta (p),$$ and consequently $\Phi_i(p) = \Theta (p)$ for every projection $p\in \CP (\mathcal{V}_{n_i})$.   
	
	Finally, we repeat the above argument with each $i\in \Gamma_1$, which produces a family of Jordan $^*$-isomorphisms $\Phi_i:  \mathcal{V}_{n_i} \to \mathcal{W}_{n_{\sigma(i)}}$. The mapping $\Phi = (\Phi_i)_{i\in \Gamma_1}$ is the desired Jordan $^*$-isomorphism. 
\end{proof}

Let $\Theta:\CP(\CA) \to \CP(\CB)$ be an order isomorphism preserving points at diametrical distance between the projection lattices of two JBW$^*$-algebras. We have seen above the difficulties of extending $\Theta$ to a Jordan $^*$-isomorphism from $\CA$ onto $\CB$. However, the next result shows that if $\Theta$ preserves distances between projections and both JBW$^*$-algebras have atomic type $I_2$ parts, the extension is possible. 

\begin{theorem}\label{thm:ord-pre-isom} Let $\CA$ and $\CB$ be JBW$^*$-algebras, and let $\Theta:\CP(\CA) \to \CP(\CB)$ be an order isomorphism which is also an isometry. Suppose additionally that the type $I_2$ part of $\CA$ is atomic. Then $\Theta$ extends to a Jordan $^*$-isomorphism $\Phi: \CA\to \CB$.
\end{theorem}

\begin{proof} Let $e_{I_2}$ and $f_{I_2}$ be the central projections in $\CA$ and $\CB$ corresponding to the type $I_2$ parts of $\CA$ and $\CB$, respectively. \Cref{prop 4 point 1}$(a)$ assures that
	$$\Theta(p) = \Theta(e_{I_2}\circ p) \vee \Theta((\one-e_{I_2})\circ p) = \Theta(e_{I_2}\circ p) + \Theta((\one-e_{I_2})\circ p) \qquad (p\in \CP(\CA)).$$
	By \Cref{prop 4 point 1}$(c)$, the map $e_{I_2} \circ p\mapsto \Theta(e_{I_2}\circ p)$ maps $\CP({\U_{e_{I_2}}(\CA)})$ onto $\CP({\U_{f_{I_2}}(\CB)})$, the mapping $(\one-e_{I_2})\circ p\mapsto \Theta((\one-e_{I_2})\circ p)$ will send $\CP({\U_{\one-e_{I_2}}(\CA)})$ onto $\CP({\U_{\one-f_{I_2}}(\CB)})$, and both of them are order isomorphisms preserving distances by the hypotheses in this theorem.  \Cref{prop 4 point 1}$(b)$ shows that the map $(\one-e_{I_2})\circ p\mapsto \Theta((\one-e_{I_2})\circ p)$ extends to
	a Jordan $^*$-isomorphism $\Phi_2: \U_{\one-e_{I_2}}(\CA)\to \U_{\one-f_{I_2}}(\CB)$. We are thus left with the case of JBW$^*$-algebras of type ${I}_2$, {whose} domain is atomic.
	
	Consider the mapping $\Theta|_{_{\CP({\U_{e_{I_2}}(\CA)})}}: \CP({\U_{e_{I_2}}(\CA)})\to \CP({\U_{f_{I_2}}(\CB)})$, where $\U_{e_{I_2}}(\CA)$ and $\U_{f_{I_2}}(\CB)$ are type ${I}_2$ JW$^*$-algebras and $\U_{e_{I_2}}(\CA) = e_{I_2}\circ \CA$ is atomic by the hypotheses. Then $e_{I_2}\circ \CA$ can be decomposed as the direct sum of a family of spin factors, that is, $e_{I_2}\circ \CA = \bigoplus_{i\in \Gamma_1}^{\ell_{\infty}} {\mathcal{V}_i}$ where each ${\mathcal{V}_i}$ is a spin factor. 
	
	Let $\one_i$ denote the unit element in ${\mathcal{V}_i}$. Clearly, minimal projections in $\Z(e_{I_2}\circ \CA)$ are mapped to minimal projections in $\Z(f_{I_2}\circ \CB)$ by $\Theta$ (cf. \Cref{prop 4 point 1}), and by the assumptions $\CP(\Z(e_{I_2}\circ \CA)) = \{ w^*\hbox{-}\sum_{i\in \Gamma_1} \sigma_{i} \one_i : i\in \Gamma_1, \sigma_i \in \{0,1\}\}$. Moreover, since $\Theta$ is an order isomorphism, $$\Theta \left(w^*\hbox{-}\sum_{i\in \Gamma_1} \sigma_{i} \one_i\right) = w^*\hbox{-}\sum_{i\in \Gamma_1} \sigma_{i} \Theta(\one_i),$$  where the $\Theta(\one_i)$'s are mutually orthogonal minimal central projections in $f_{I_2}\circ \CB$. Therefore $f_{I_2}\circ \CB$ must be also atomic (and of type $I_2$). Then the mapping  $\Theta|_{_{\CP({\U_{e_{I_2}}(\CA)})}}: \CP({\U_{e_{I_2}}(\CA)})\to \CP({\U_{f_{I_2}}(\CB)})$ is an isometric order isomorphism. So, by the final statement in \Cref{p atomic type I2}{,} $\Theta|_{\CP (\CA_{2})}$ admits a unique extension to a Jordan $^*$-isomorphism $\Phi_1: e_{I_2}\circ \CA\to f_{I_2}\circ \CB$. To conclude the proof it suffices to take $\Phi = \Phi_1 +\Phi_2$.  
\end{proof}

To {culminate} this section we gather some characterizations of algebraic properties in terms of the relation ``being at diametrical distance''. Observe first that for each subset $E\subseteq \RS_{\CA^+}$ we have 
\begin{align*}
	\sphe(E) & = \big\{ b\in \RS_{\CA^+} :  b \hbox{ is at diametrical distance from every } e\in E \big\}.
\end{align*}

\begin{corollary}\label{c: rem:inver} Let $\CA$ be a JBW$^*$-algebra. The following statements hold:
	\begin{enumerate}[$(a)$] \item $\CP(\CA) \setminus \{0\} = \left\{a\in \RS_{\CA^+}: \sphe(\sphe(a)) = \{a\} \right\}$.
		\item The unit element $\one \in {\CA^+}$ can be characterized as the unique element $a\in \RS_{\CA^+}$ satisfying the following properties:
		\begin{enumerate}[$(1)$]
			\item $\{a\} = \sphe(\sphe(a))$;
			\item $\CP(\CA) \setminus \{0,a\}\subseteq \sphe(a)$ $($equivalently, $\| a - q\| = 1$ for every $q\in \CP(\CA)\setminus \{a\})$;
			\item $\sphe(b)\cup \sphe(a)\neq \RS_{\CA^+}$, for every $b\in \RS_{\CA^+}$.
		\end{enumerate}
		\item $\label{eqt:inver}
		\RS_{\CA^+}^{-1} = \{c\in \RS_{\CA^+}: \|\one-c\|<1 \} = \RS_{\CA^+} \setminus \sphe(\one)$. 
		\item $q\leq p$ in $\CP(\CA) \setminus \{0\}$ if, and only if, $\sphe (p)\cap  \left(\RS_{\CA^+} \setminus \sphe(\one)\right) \subseteq \sphe (q).$
		\item Given $q, p$ in $\CP(\CA)\setminus \{0\}$ we have $p$ and $q$ are orthogonal if, and only if, for every $r,s$ in $\CP(\CA)\setminus \{0\}$ with $r\leq p$ and $s\leq q$ we have $r\in \sphe(s)$ and $s\in \sphe(r)$. 
	\end{enumerate}  
\end{corollary}

\begin{proof} The statement in $(a)$ is a rephrasing of \Cref{prop:lem:metr-discr-proj}$(b)$. The just quoted proposition and \Cref{Lemma 2point2} give the statement in $(b)$. The identity in $(c)$ is clear. The equivalence in $(d)$ follows from statement $(c)$ above and \Cref{lem:describ-ord-proj}. Finally \Cref{lem:describ-disj-proj} proves $(e)$.
\end{proof}

The next corollary is a straightforward consequence of the previous result since all types of elements and relationships among them described in the previous corollary are given by the relationship ``being at diametrical distance''.  

\begin{corollary}\label{c:lem:met-pre-proj} Let $\Delta : \RS_{\CA^+} \to \RS_{\CB^+}$ be a bijection between the positive unit spheres of two JBW$\,^*$-algebras. Suppose that $\Delta$ preserves points at diametrical distance, i.e., $$\|a-b\| = 1 \Leftrightarrow \|\Delta(a)-\Delta(b)\| = 1.$$ This property is also known as the strong distance-one preserving property (SDOPP) in the literature (see, for example, \cite{Ding2006}). Then the following statements hold:
	\begin{enumerate}[$(a)$]
		\item $\Delta(\CP (\CA)\setminus \{0\}) = \CP(\CB)\setminus \{0\}$.
		\item $\Delta(\one) = \one$ and $\Delta\big(\RS_{\CA^+}^{-1}\big) = \RS_{\CB^+}^{-1}$.
		\item If we extend $\Delta|_{\CP(\CA)\setminus \{0\}}$ to a map $\bar{\Delta} : \CP(\CA)\to \CP(\CB)$ by setting $\bar{\Delta} (0) := 0$, then $\bar {\Delta}$ is an order isomorphism preserving points at diametrical distance. Furthermore, if $\Delta$ is an isometry, {the mapping} $\bar {\Delta}$ enjoys the same property. 
	\end{enumerate}
\end{corollary}

\section{Surjective isometries between positive unit spheres of JBW$^*$-algebras}\label{Sec:isometries between positive spheres}

This section is fully devoted to the study of surjective isometries between the positive unit spheres of two JBW$^*$-algebras.   

\begin{proposition}\label{prop:tol-on-S+=>Jord} Let $\CA$ and $\CB$ be JBW$\,^*$-algebras. Assume that the type $I_2$ part of $\CA$ is atomic. Suppose that there {exists} a bijection $\Delta: \RS_{\CA^+} \to \RS_{\CB^+}$ preserving points at distances $1$ and $\frac{\sqrt{2}}{2}$ in both directions. Then the following statements hold:
	\begin{enumerate}[$(a)$]
		\item $\CA$ is Jordan $^*$-isomorphic to $\CB$.
		\item If $\CA$ has no type {$I_2$} part, then $\Delta|_{\CP(\CA)\setminus \{0\}}$ extends to a Jordan $^*$-isomorphism from $\CA$ onto $\CB$.
	\end{enumerate}
\end{proposition}

\begin{proof} \Cref{c:lem:met-pre-proj}$(c)$ assures that $\Delta|_{\CP(\CA)\setminus \{0\}}$ extends, in a canonical way, to a map $\bar{\Delta} : \CP(\CA)\to \CP(\CB)$ which is an order isomorphism preserving points at diametrical distance. \Cref{prop 4 point 1}$(c)$ proves that $\bar{\Delta}$ maps the projections in the type $I_2$ part (respectively, the orthogonal complement of the type $I_2$ part) of $\CA$ onto the projections in the type $I_2$ part (respectively, the orthogonal complement of the type $I_2$ part) of $\CB$. \Cref{prop 4 point 1}$(b)$ implies that the orthogonal complement of the type $I_2$ part of $\CA$ is Jordan $^*$-isomorphic to  the orthogonal complement of the type $I_2$ part of {$\CB$}. As in the proof of \Cref{thm:ord-pre-isom}, the type $I_2$ part of $\CB$ is also atomic. The first statement of \Cref{p atomic type I2} now shows that the type $I_2$ parts of $\CA$ and $\CB$ are Jordan $^*$-isomorphic. This concludes the proof of $(a)$.
	
	The statement in $(b)$ follows from \Cref{c:lem:met-pre-proj}$(b)$.
\end{proof}

The next lemma is an extension of \cite[Lemma 2.3(c)]{LN} to JBW$^*$-algebras. The conclusion of this result, together with \Cref{c: rem:inver}$(c)$ and \Cref{prop:lem:metr-discr-proj}, shows how to determine whether a projection $q$ is dominated by an invertible element $a\in \RS_{\CA^+}^{-1}$ in terms of the relation ``being at diametrical distance''.

\begin{proposition}\label{l: second lemma p leq a in section 5} Let $\CA$ be {a JBW$\,^*$-algebra}. Suppose that $q\in \CP(\CA)\setminus \{0,\one\}$ and $a\in \RS_{\CA^+}^{-1}$. Then the following statements are equivalent:
	\begin{enumerate}[$(a)$]
		\item $q\leq a$.
		\item $U_q (a) =q$ is the unit of $U_{q}(\CA)$, equivalently, $a = q + U_{\one-q} (a)\in q + \RB_{U_{\one-q}(\CA)^+}^{-1}$. 
		\item $U_q (a) \in \RS_{\CA^+}$ and satisfies the following properties:
		\begin{enumerate}[$(1)$]
			\item $\{U_q (a)\} = \spheAq(\spheAq(U_q (a)))$;
			\item $\CP(U_q(\CA)) \setminus \{0,U_q (a)\}\subseteq \spheAq(U_q (a))$ $($equivalently, $\| U_q (a) - r\| = 1$ for every $r\in \CP(U_q(\CA))\setminus \{U_q (a)\})$;
			\item $\spheAq(b)\cup \spheAq(U_q (a))\neq \RS_{U_q(\CA)^+}$, for every $b\in \RS_{U_q(\CA)^+}$.
		\end{enumerate}
		\item $\|\one-r-a\| = 1$ $($i.e., $\one-r\in \sphe(a) )$ for every $r\in \CP(\CA)\setminus \{0\}$ with $r\leq q$. 
	\end{enumerate}	
\end{proposition}

\begin{proof} $(a)\Rightarrow (b)$ Suppose $\one\geq a\geq q$. Since the mapping $U_q$ is positive (cf. \cite[Proposition 3.3.6]{HOS}) it follows that $q\leq U_q (a) \leq U_q (\one ) = q$. The equivalence with the second statement follows from \Eqref{eq FR Peirce 2 plus orthogonal}.
	
	The equivalence $(b)\Leftrightarrow (c)$ is a consequence of \Cref{c: rem:inver}$(b)$.   
	
	$(b)\Rightarrow (d)$ If $a= q + U_{\one-q} (a)$, by orthogonality, for each $0\neq r \leq q$ we have $\|\one -r - a\| = \|- r +(\one - q) - U_{\one-q} (a)\|=\max\{\|r\|,\|(\one - q) - U_{\one-q} (a)\|\}= 1$.	
	
	$(d)\Rightarrow (a)$ The element $U_q (a)$ is positive in $U_q(\CA)$. The invertibility of $a$ implies that $U_q (a)\neq 0$. Let $\mathcal{C}$ stand for the JBW$^*$-subalgebra of $\CA$ generated by $q$ and $U_q(a)$. It is known that $\mathcal{C}$ is a commutative von Neumann algebra. If $U_q (a)\lneq q$, by standard functional calculus, there exists a projection $0\neq r\leq q$ such that $\| r\circ U_q (a)\| = \|U_r U_q(a)\|<1$.
	
	By the assumptions $\|\one-r -a\| =1$. Having in mind that $a$ is invertible, by applying \Cref{lem:peak-and-zero}$(a)$ and $(b)$ we deduce the existence of $\omega\in \ps(\CA)$ satisfying $\omega (a)=1$ and $\omega (\one-r)=0$. \Cref{lem:peak-and-zero}$(c)$ assures that $\omega (x) = \omega (U_r(x))$ for all $x\in \CA$. Therefore, $$1=\omega (a) = \omega (U_r(a)) = \omega (U_rU_q(a))\leq \|U_r U_q (a)\| <1,$$ which is impossible.   
\end{proof}

We are now in a position to establish a positive answer to the so-called Tingley's problem for positive spheres of JBW$^*$-algebras whose type $I_2$ part is atomic.

\begin{theorem}\label{thm:met-pre-proj} Let $\CA$ and $\CB$ be JBW$\,^*$-algebras such that the type $I_2$ part of $\CA$ is atomic. Let $\Delta:\RS_{\CA^+} \to \RS_{\CB^+}$ be a surjective isometry. Then $\Delta$ extends $($uniquely$)$ to a Jordan $^*$-isomorphism from $\CA$ onto $\CB$.
\end{theorem}

\begin{proof}
	By \Cref{c:lem:met-pre-proj}$(a)$ and $(c)$ and \Cref{thm:ord-pre-isom}, the bijection  $$\Delta|_{\CP(\CA)\setminus \{0\}}:\CP(\CA)\setminus \{0\} \to \CP(\CB)\setminus \{0\}$$ extends to a Jordan $^*$-isomorphism $\Phi: \CA \to \CB$. It remains to show that $\Phi|_{\RS_{\CA^+}} = \Delta$.
	
	\Cref{c:lem:met-pre-proj}$(b)$ and \Cref{l: second lemma p leq a in section 5} assure that the following statements hold:
	\begin{equation}\label{eqt:pres-order}
		\begin{aligned}
			&\Delta(\one) =\one, \hbox{ and for all } p\in \CP({A})\setminus \{0,\one\}, \hbox{ and all } a\in \RS_{\CA^+}^{-1} \\	
			&\hbox{ we have } p\leq a \text{ if, and only if, } \Delta(p)\leq \Delta(a),
		\end{aligned}
	\end{equation} essentially because $\Delta(\one-r) = \Phi(\one-r) = \one -\Delta(r)$ and $\Delta(r) = \Phi(r)\leq \Delta(p)$ when $r\in \CP({A})\setminus \{0\}$ with $r\leq p$.
	
	Consider $p_0\in \CP(\CA)\setminus \{0,\one\}$, and set $q_0:= \Delta(p_0)$. Following the standard notation associated with the Peirce decomposition, we write $$\CA_0 (p_0):= \U_{\one-p_0}(\CA)\quad \text{and} \quad \CB_{0}(q_0):= \U_{\one-q_0}(\CB).$$ Let us note that $\CA_0 (p_0)$ and $\CB_{0}(q_0)$ are unital JBW$^*$-algebras with units $\one-p_0$ and $\one-q_0$, respectively.  
	
	We shall next prove the existence of a Jordan $^*$-isomorphism $\Phi^{p_0}: \CA_0 (p_0) \to \CB_0 (q_0)$  satisfying
	\begin{align}\label{eqt:Lambda0a}
		\Delta(p_0 + a) &= q_0 +\Phi^{p_0}(a)\quad \big(a\in \RB_{{\CA_0 (p_0)}^+}\big).
	\end{align}
	Indeed, since $p_0+ \RB_{{\CA_0 (p_0)}^+}^{-1} = \left\{ a\in \RS_{\CA^+}^{-1} : a\geq p_0\right\}$ and $q_0+ \RB_{{\CB_0 (q_0)}^+}^{-1} = \left\{ b\in \RS_{\CB^+}^{-1} : b\geq q_0\right\}$, it follows from \Cref{l: second lemma p leq a in section 5} that 
	$$
	\Delta\big(p_0+ \RB_{{\CA_0 (p_0)}^+}^{-1}\big) = q_0 + \RB_{{\CB_0 (q_0)}^+}^{-1}.
	$$
	Since $\Delta$ is an isometry and the sets $\RB_{{\CA_0 (p_0)}^+}^{-1}$ and $\RB_{{\CB_0 (q_0)}^+}^{-1}$ are {norm-dense} in $\RB_{{\CA_0 (p_0)}^+}$ and $\RB_{{\CB_0 (q_0)}^+}$, respectively, we have
	$$
	\Delta\big(p_0+\RB_{{{\CA_0 (p_0)}}^+}\big) = q_0 +\RB_{{\CB_0 (q_0)}^+}.
	$$
	Up to composing with the translations by the vectors $p_0$ and $q_0$, we can find an isometric bijection $\Delta^{p_0}: \RB_{{\CA_0 (p_0)}^+} \longrightarrow \RB_{{\CB_0 (q_0)}^+}$ satisfying $\Delta(p_0 + a) = q_0 +\Delta^{p_0}(a),$ for all $a\in \RB_{{\CA_0 (p_0)}^+}$. Having in mind that $\RB_{{\CA_0 (p_0)}^+}$ and $\RB_{{\CB_0 (q_0)}^+}$ are closed convex subsets with non-empty interior, we know from \cite[Theorem 5]{mankiewicz1972extension} that $\Delta^{p_0}$ extends to a bijective affine isometry, denoted by $\Phi^{p_0}$, from $\CA_0 (p_0)$ onto $\CB_0 (q_0)$. Since $\Delta^{p_0}(0) = 0$ and $\Delta^{p_0} (\one-p_0) =  \one-q_0$, the mapping $\Phi^{p_0}$ is necessarily linear ($\Phi^{p_0} (0) =0$) and unital ($\Phi^{p_0} (\one-p_0)=\one - q_0$). Therefore $\Phi^{p_0}: \CA_0 (p_0)\to \CB_0 (q_0)$ is a unital and surjective real linear isometry between two JBW$^*$-algebras.  Corollary 3.4 in \cite{FerMarPe2004} assures that $\Phi^{p_0}$ is a Jordan $^*$-isomorphism.
	
	Next, we claim that $\Phi^{p_0}$ satisfies the following property:
	\begin{equation}\label{eqt:Jp0=Lambda}
		\Phi^{p_0}(e) = \Delta(e) = \Phi(e) \quad (e\in \CP({\CA_0 (p_0)})\setminus \{0\}).
	\end{equation} In fact, consider $e\in \CP({\CA_0 (p_0)})\setminus \{0\}$. The equality $\Delta(e) = \Phi(e)$ comes from the definition of $\Phi$. On the other hand, the identity in \Eqref{eqt:Lambda0a} gives $$q_0 + \Phi^{p_0}(e) = \Delta(p_0 +e) = \Phi(p_0+e) = q_0 + \Phi(e),$$
	and hence $\Phi^{p_0}(e) =  \Phi(e)$.

	Finally, consider $a\in \RS_{\CA^+}$ and $\varepsilon > 0$. Working on the JBW$^*$-subalgebra of $\CA$ generated by $a$ and $\one$, which is a commutative von Neumann algebra, one can find $\alpha_1,\dots, \alpha_n$ in $(0,1),$ as well as pairwise orthogonal projections $p_0,p_1,\dots,p_n\in \CP(\CA)\setminus \{0\}$ such that $\left\|a - \left(p_0 + \sum_{i=1}^n \alpha_i p_i\right)\right\| < \varepsilon$.
	As both $\Delta$ and $\Phi$ are isometries, we have
	$$\left\|\Delta(a) - \Delta\left(p_0 + {\sum}_{i=1}^n \alpha_i p_i\right)\right\| < \varepsilon,
	\quad \text{and}\quad \left\|\Phi(a) - \Phi\left(p_0 + {\sum}_{i=1}^n \alpha_i p_i\right)\right\| < \varepsilon.$$
	It is clear that ${\sum}_{i=1}^n \alpha_i p_i \in \RB_{{\CA_0 (p_0)}^+}$.
	Thus, {the} identities in \Eqref{eqt:Lambda0a} and \Eqref{eqt:Jp0=Lambda}, as well as the definition of $\Phi$ imply that
	\begin{align*}
		\Delta\left(p_0 + {\sum}_{i=1}^n \alpha_i p_i\right)
		& =  \Delta(p_0) + \Phi^{p_0}\left({\sum}_{i=1}^n \alpha_i p_i\right)
		= q_0 + {\sum}_{i=1}^n \alpha_i \Phi(p_i)\\
		& = \Phi(p_0) + {\sum}_{i=1}^n \alpha_i \Phi(p_i)
		= \Phi\left(p_0 + {\sum}_{i=1}^n \alpha_i p_i\right).
	\end{align*} The arbitrariness of $\varepsilon>0$ gives $\Delta(a) = \Phi(a)$, as required.  
\end{proof}

\medskip

\section*{Acknowledgements}

\noindent A.M. Peralta supported by grant PID2021-122126NB-C31 funded by MICIU/AEI/\linebreak 10.13039/501100011033 and by ERDF/EU, by Junta de Andalucía grant FQM375, IMAG--Mar{\'i}a de Maeztu grant CEX2020-001105-M/AEI/10.13039/501100011033 and (MOST) Ministry of Science and Technology of China grant G2023125007L. \smallskip

\noindent P. Saavedra supported by a ``Formaci\'{o}n de Profesorado Universitario (FPU)'' grant from the  Spanish Ministry of Science, Innovation, and Universities (Grant No. FPU23/00747).\smallskip

\noindent We thank C.-W. Leung, C.-K. Ng, N.-C. Wong for fruitful discussions and comments on the manuscript.\medskip

\smallskip\smallskip

\noindent\textbf{Data Availability} Statement Data sharing is not applicable to this article as no datasets were generated or analysed during the preparation of the paper.\smallskip\smallskip

\noindent\textbf{Declarations} 
\smallskip\smallskip

\noindent\textbf{Conflict of interest} The authors declare that they have no conflict of interest.

\end{document}